\documentclass[12pt,authoryear]{elsarticle}
\usepackage{amssymb}
\usepackage{amsmath}
\usepackage{bm}
\usepackage{epsfig,subfigure}
\usepackage[norelsize,lined,algoruled]{algorithm2e}
\usepackage{color}
\usepackage{float}
\usepackage{hyperref}
\usepackage{textgreek}
\usepackage[title]{appendix}
\usepackage{url}

\SetKwInput{KwInput}{Input}
\SetAlCapSkip{1em}
\SetKwFunction{defined}{defined}
\SetKwFunction{notdefined}{not defined}
\SetKwFunction{mean}{mean}
\SetKwProg{Fn}{Function}{}{}
\SetKwFunction{FnGetLowerBound}{GetAlphaLowerBound}%
\SetKwFunction{FnGetAlpha}{GetAlphaEstimate}%

\newtheorem{defn}{Definition}
\newtheorem{remark}{Remark}
\newcommand{\dotstar}{\cdot \ast}

\newcommand{\bfa}{\mathbf{a}}
\newcommand{\bfb}{\mathbf{b}}
\newcommand{\bfc}{\mathbf{c}}
\newcommand{\bcext}{\bfc^{\mathrm{ext}}}
\newcommand{\barc}{\bar{\bfc}}
\newcommand{\barcext}{\barc^{\mathrm{ext}}}
\newcommand{\bfd}{\mathbf d}
\newcommand{\dx}{\Delta_x}
\newcommand{\dy}{\Delta_y}
\newcommand{\dz}{\Delta_z}

\newcommand{\Gr}{G^{(r)}}

\newcommand{\barGr}{\bar{G}^{(r)}}
\newcommand{\barg}{\bar{g}}
\newcommand{\bargamma}{\bar{\gamma}}

\newcommand{\bfm}{\mathbf m}
\newcommand{\nbr}{n_r}
\newcommand{\nbx}{n_x}
\newcommand{\nby}{n_y}
\newcommand{\nbz}{n_z}
\newcommand{\nsx}{s_x}
\newcommand{\nsy}{s_y}
\newcommand{\pady}{p_y}
\newcommand{\padx}{p_x}
\newcommand{\padxl}{p_{x_{\textrm{L}}}}
\newcommand{\padxr}{p_{x_{\textrm{R}}}}
\newcommand{\padyl}{p_{y_{\textrm{L}}}}
\newcommand{\padyr}{p_{y_{\textrm{R}}}}
\newcommand{\param}{\mathtt{prob\_params}}
\newcommand{\bfr}{\mathbf{r}}
\newcommand{\brext}{\bfr^{\mathrm{ext}}}
\newcommand{\barr}{\bar{\bfr}}
\newcommand{\barrext}{\barr^{\mathrm{ext}}}
\newcommand{\Rmn}[2]{\mathcal{R}^{#1\times #2}}
\newcommand{\Rm}[1]{\mathcal{R}^{#1}}
\newcommand{\That}{\mathrm{That}}
\newcommand{\Thatforward}{\mathrm{That.forward}}
\newcommand{\Thattranspose}{\mathrm{That.transpose}}
\newcommand{\Tcirc}{T^{\mathrm{circ}}}

\newcommand{\bfu}{\mathbf{u}}

\newcommand{\bfv}{\mathbf{v}}
\newcommand{\bfw}{\mathbf{w}}
\newcommand{\bfx}{\mathbf{x}}
\newcommand{\bfy}{\mathbf{y}}
\newcommand{\txtrho}{\text{\textRho}}
\newcommand{\txtups}{\text{\textUpsilon}}
\newcommand{\txtchi}{\text{\textChi}}
\newcommand{\txtzet}{\text{\textZeta}}

\newcommand{\bvec}{\texttt{vec}}
\newcommand{\barray}{\texttt{array}}
\newcommand{\flipud}{\texttt{flipud}}
\newcommand{\fliplr}{\texttt{fliplr}}
\newcommand{\fftt}{\texttt{fft2}}
\newcommand{\ifftt}{\texttt{ifft2}}
\newcommand{\fftw}{\texttt{fftw}}
\newcommand{\dwisdom}{\texttt{dwisdom}}
\newcommand{\real}{\texttt{real}}

\newcommand{\reshape}{\texttt{reshape}}
\newcommand{\bttb}{\texttt{BTTB}}
\newcommand{\multbttb}{\texttt{mult\_BTTB}}
\newcommand{\symbttb}{\texttt{sym\_BTTB}}
\newcommand{\symbttbfft}{\texttt{sym\_BTTBFFT}}
\newcommand{\bttbfft}{\texttt{BTTBFFT}}
\newcommand{\gravity}{\texttt{gravity}}
\newcommand{\magnetic}{\texttt{magnetic}}
\newcommand{\Ggrav}{G_{\gravity}}
\newcommand{\Gmag}{G_{\magnetic}}
\newcommand{\Tgrav}{T_{\gravity}}
\newcommand{\Tmag}{T_{\magnetic}}

\begin{document}
\begin{frontmatter}

\title{A Tutorial and Open Source Software for the Efficient Evaluation of Gravity and Magnetic Kernels}
\author[1]{Jarom D. Hogue}
\ead{jdhogue@asu.edu}
\author[2]{Rosemary Anne Renaut\corref{cor1}}
\ead{renaut@asu.edu}
\author[3]{Saeed Vatankhah }
\ead{svatan@ut.ac.ir}

\address[1]{School of Mathematical and Statistical Sciences, Arizona State University, Tempe, AZ, USA}
\address[2]{School of Mathematical and Statistical Sciences, Arizona State University, Tempe, AZ, USA}
\address[3]{Institute of Geophysics, University of Tehran, Tehran, Iran}
\cortext[cor1]{Corresponding author}

\begin{abstract}
Fast computation of three-dimensional gravity and magnetic forward models is considered. When the measurement data is assumed to be obtained on a uniform grid which is staggered with respect to the discretization of the parameter volume, the resulting kernel sensitivity matrices exhibit block-Toeplitz Toeplitz-block (BTTB) structure. These matrices are symmetric for the gravity problem but non-symmetric for the magnetic problem. In each case, the structure facilitates fast forward computation using two-dimensional  fast Fourier transforms. The construction of the kernel matrices and the application of the transform for fast forward multiplication,  for each problem, is carefully described. But, for purposes of comparison with the transform approach, the generation of the unique entries that define a given kernel matrix is also explained. It is also demonstrated how the matrices, and hence transforms, are adjusted when padding around the volume domain is introduced.  The transform algorithms for fast forward matrix multiplication with the sensitivity matrix and its transpose, without the direct construction of the relevant matrices,  are presented. Numerical experiments demonstrate the significant reduction in computation time that is achieved using the transform implementation.  Moreover, it becomes feasible, both in terms of reduced memory requirements and computational time,  to implement the transform algorithms for large three-dimensional volumes. All presented algorithms, including with variable padding,  are coded for optimal memory, storage and computation as an open source MATLAB code which can be adapted for any convolution kernel which generates a BTTB matrix, whether or not it is symmetric. This work, therefore, provides a general  tool for the efficient simulation of gravity and magnetic field data, as well as any formulation which admits a sensitivity matrix with the required structure. 

 \end{abstract}

\begin{keyword}
Forward modeling \sep fast Fourier Transform\sep  Gravity\sep Magnetic
\end{keyword}
\date{\today}

\end{frontmatter}
\newpage

\section{Introduction}\label{sec:introduction}

Fast computation of geophysics kernel models has been considered by a number of authors, including calculation within the Fourier domain as in \cite{Li2018,Pilkington:97,ZHAO2018294}, and through discretization of the operator and calculation in the spatial domain as in \cite{ChenLiu:18,ZhangWong:15}. \cite{Pilkington:97} introduced the use of the  Fast Fourier Transform (FFT) for combining the evaluation of the magnetic kernel in the Fourier domain with the conjugate gradient method for solving the inverse problem to determine magnetic susceptibility from  measured magnetic field data.  \cite{Li2018}  considered the use of the Gauss FFT for fast forward modeling of the magnetic kernel on an undulated surface, combined with spline interpolation of the surface data. This work focused on the implementation of the model in the wave number domain and only applied the method for forward modeling. The Gauss FFT was also used by \cite{ZHAO2018294} for the development of a high accuracy forward modeling approach  for the gravity kernel. 

\cite{bruun2007}, and subsequently, \cite{ZhangWong:15},  introduced the use of the Block-Toeplitz-Toeplitz-Block (BTTB) structure of the modeling sensitivity matrix for fast three-dimensional inversion of three-dimensional gravity and magnetic data. For a matrix with BTTB structure, it is possible to embed the information within a matrix of Block-Circulant Circulant-Block (BCCB) structure that facilitates fast forward multiplication using a two-dimensional FFT (2DFFT), see e.g. \cite{ChanFuJin:2007,Vogel:2002}.  For three-dimensional modeling, \cite{ZhangWong:15}  exploited the two-dimensional multi-layer structure of the kernel, that provides BTTB structure for each layer of the domain, and performed the  inverse operation iteratively over all layers of the domain. The technique is flexible to depth layers of variable heights, and permits the inclusion of smoothness stabilizers in the inversion, for each layer of the domain. Moreover, \cite{ZhangWong:15} adopt the preconditioning of the BTTB matrix using  optimal preconditioning operators as presented in \cite{ChanFuJin:2007} for implementing efficient and effective solvers for the inversion. 

A fast forward modeling of the gravity field was  developed by  \cite{ChenLiu:18}, using the three-dimensional modeling of the gravity kernel as given in  e.g.  \cite{BoCh:2001,haaz1953}. Their work extends the techniques of \cite{ZhangWong:15} for taking advantage of the BTTB structure of the sensitivity matrix associated with a single layer of the gravity kernel, but offers greater improvements in the implementation of the forward kernel through the presentation of an optimized calculation of the kernel entries in this matrix. This arises due to the observation that the uniform placement of the measurement stations, in relation to the coordinate grid for the unknown densities, yields significant savings in computation and memory. When the stations are on a grid that is uniformly staggered with respect to the coordinate grid on the top surface of the coordinate domain,  redundant operations in the calculation of the  the sensitivity matrix can be eliminated.

Here, we present a careful derivation of the forward operators for kernels that are spatially invariant in all dimensions, and are thus convolution operators. With uniform placement of measurement stations relative to the coordinate domain, namely on staggered grids in the $x$ and $y$ dimensions, the resulting discretizations of the first kind Fredolm integral operators, \cite{Zhd:2002}, yield sensitivity matrices that exhibit BTTB  structure. We distinguish between operators that yield symmetric BTTB (symBTTB) matrices and those that are BTTB but are not symmetric. This depends on the kernel operator, the gravity kernel, \cite{BoCh:2001}, is symmetric in the distances, but the magnetic kernel, \cite{RaoBabu:91} is not. Thus, the generation of the sensitivity matrix requires further analysis for efficient computation in the case of the magnetic kernel as compared to the gravity kernel. Here we demonstrate that it is feasible to develop an optimized calculation of the entries of the magnetic kernel matrix, in a manner similar to that used for the gravity case, but due to lack of symmetry the memory and computation requirements are increased. Still, remarkable savings in generating the matrix are achieved.   We note further, that while it is not possible to take advantage of BTTB structure when the stations are not on a uniform grid, the calculation of the underlying kernel matrices can  still be optimized for arbitrary stations locations, by reuse of common vectors and arrays for each depth layer of the coordinate volume.

We consider the general case for matrices with BTTB structure and then relate the discussion to the specific geophysics gravity and magnetic kernels. Although these topics have been discussed in the literature, the analysis for the magnetic kernel is new. Moreover, this work provides both the careful derivation of the matrices that arise, and the associated implementation of the FFT for fast computation of forward multiplications with the matrix and its transpose. The intent is to provide a user-friendly environment for the development and validation of kernel operators that yield the underlying BTTB structures which facilitate the fast computation of both  operator and forward operations.

\textit{Overview of main scientific contributions.}  Our approach implements and extends the BTTB algorithm for the forward modeling with the magnetic kernel, and for the inclusion of padding around the domain.  Specifically, our main contributions are as follows.
(i) We present a detailed derivation of the implementation of the algorithm presented in \cite{ChenLiu:18} for the forward modeling of the gravity problem; (ii) The algorithm is extended to include domain padding in $x$ and $y$ directions, that is not necessarily symmetric with respect to the domain; (iii) We demonstrate the use of the 2DFFT for forward multiplication using the transpose matrix, as is required for the solution of the associated inverse problem; (iv) The algorithm, with and without padding, is further extended for matrices that are not symmetric, as occurs for the magnetic kernel; (v) Efficient derivation of the underlying operators to be used without the 2DFFT is also provided, so as to facilitate a realistic performance comparison between the use of the FFT for the forward multiplication and a direct forward multiplication without the use of the FFT;  (vi) All components of the discussion are coded for optimal memory, storage and computation as an open source MATLAB code which can be adapted for any convolution kernel which generates a BTTB matrix. The user need only provide the algorithm for the computation of the original matrix components.  

In summary, this document provides a detailed tutorial on the efficient generation of sensitivity matrices with BTTB structure and their efficient forward multiplication using the 2DFFT.  The approach depends on the underlying convolution kernel that describes the model. This work, therefore, provides a general  tool for the efficient simulation of gravity and magnetic field data, as well as any formulation which admits a model matrix with BTTB structure
. The algorithm is open source and available at  \url{https://github.com/renautra/FastBTTB}, along with a full description of the algorithm implementation and example simulations at \url{https://math.la.asu.edu/~rosie/research/bttb.html}. 

The paper is organized as follows. In Section~\ref{sec:forwardmodel} we present the general kernel-based forward model
, and specifically in  Section~\ref{sec:invariantkernel} for convolutional kernels as seen for gravity and magnetic potential fields.  We demonstrate  in Section~\ref{sec:stationplacement} how the placement of the measurement stations as uniformly staggered with respect to the coordinate domain yields a distance vector for distances from coordinates to stations that is efficiently stored as a one-dimensional instead of two-dimensional vector. This applies also for the case of the introduction of padding around the domain, Section~\ref{sec:stationplacementpadding}.  We then show in Section~\ref{sec:kernels} how operators that are spatially invariant yield matrix operators with BTTB structure, for the symmetric case in Sections~\ref{sec:symBTTBstructure}-\ref{sec:symBTTBstructurepadded} and the unsymmetric case in Sections~\ref{sec:BTTBstructure}-\ref{sec:BTTBstructurepadded}, where in each case we first show the case without padding and then give an example of the development of the domains with padding.  In each case we explicitly explain how the relevant entries of the matrices are calculated. In Section~\ref{sec:circulant} we show how these entries are built into the formulation that facilitates the use of the 2DFFT, following the discussion of \cite{Vogel:2002}. Specific examples are given in Section~\ref{sec:specifickernels} for the efficient derivation of the entries in the operators for gravity and magnetic kernels, following \cite{ChenLiu:18} and \cite{RaoBabu:91}, in Sections~\ref{sec:gravitykernel} and \ref{sec:magnetickernel}, respectively. We demonstrate the improved efficiency of forward operations with the matrix and its transpose for these kernels, for domains of increasing size. The presented numerical results in Section~\ref{sec:numerics} validate that the given algorithms are efficient and facilitate forward modeling for problems that are   significantly larger as compared to the case when the BTTB structure is not utilized for fast computation with the 2DFFT.  Software availability is discussed in Section~\ref{sec:code} and conclusions with topics for future work are discussed in Section~\ref{sec:conclusions}. The  adopted notation  and  algorithms are presented in  Appendices \ref{App:A} and \ref{App:B}, respectively. 

\section{Forward Modeling}\label{sec:forwardmodel}
\begin{figure}
 \begin{center}
 \includegraphics[width=.75\textwidth]{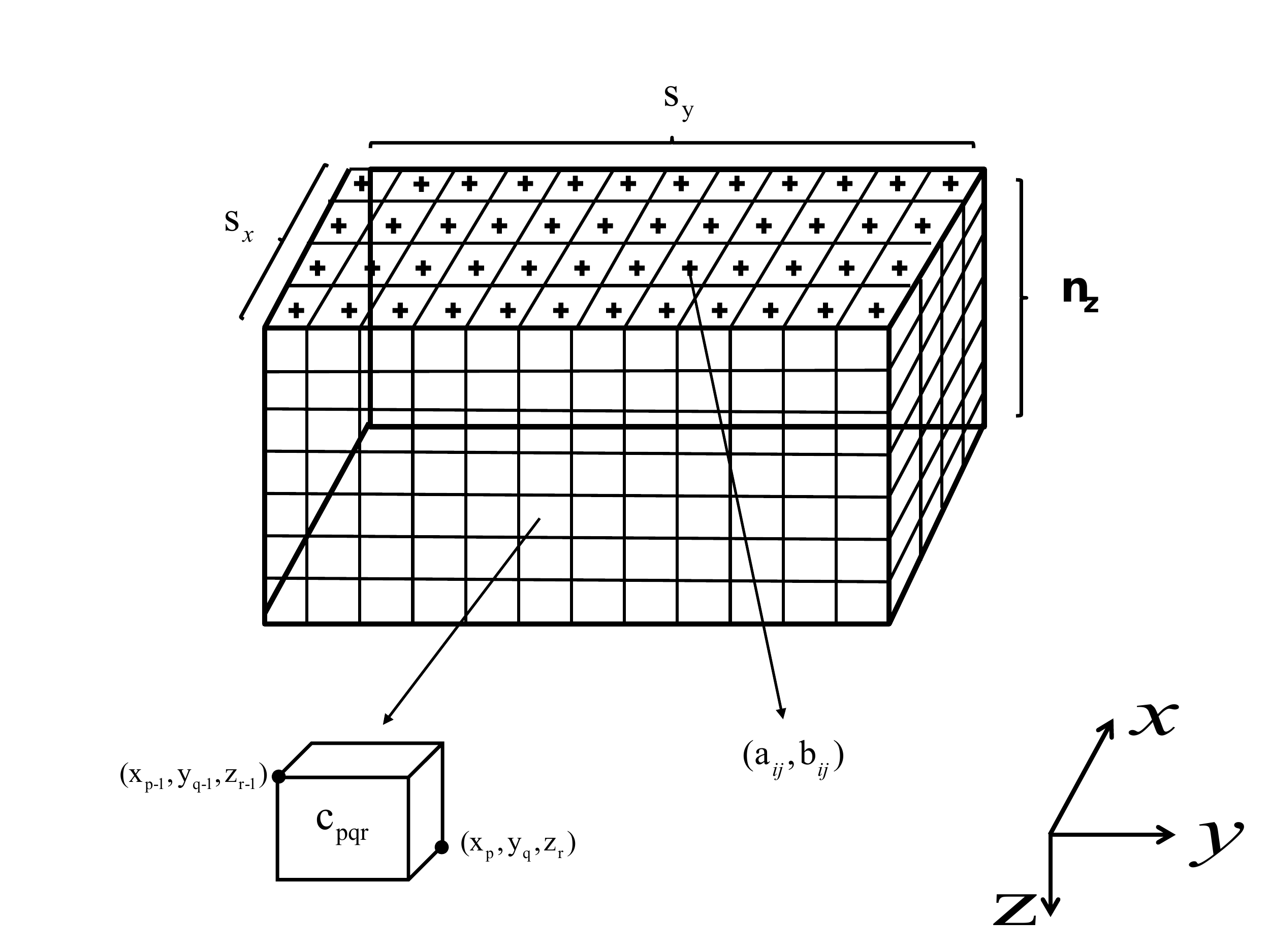}
 \caption{The configuration of prism $pqr$ in the volume relative to a station on the surface at location $s_{ij}=(a_{ij}, b_{ij})$. This shows that the stations are located on the surface of the physical domain. \label{figure1}}
 \end{center}  \end{figure}

We consider a forward model described by the Fredholm integral equation of the first kind
\begin{equation}\label{continuousforwardmodel}
d(a,b,c)=\int \int \int h(a,b,c,x,y,z) \zeta(x,y,z) dx \,dy\, dz,
\end{equation}
for which discretization leads to the forward model $\bfd=G\bfm$, with sensitivity matrix $G$, and $\bfd$ and $\bfm$ the discretizations of $d$ and $\zeta$, respectively. 
We suppose that data measurements for $d(a,b,c)$, on the surface with $c=0$,  are made at, not necessarily uniformly-spaced, station locations denoted by 
\begin{equation}\label{stationsij}
s_{ij}=(a_{ij},b_{ij}), \,\, 1\le i \le \nsx, \,\,1\le j \le \nsy.
\end{equation}
The total number of stations at the surface is $m=\nsx\nsy$.   The volume domain, without padding, is discretized into $n=\nsx\nsy\nbz$ uniform prisms, $c_{pqr}$,  with  coordinates\footnote{Note that there are, for example in the $x-$dimension, $\nsx$ blocks and hence $\nsx+1$ coordinates describing these blocks.}
\begin{equation}\label{prismspqr}
\begin{array}{ccc}
x_{p-1}=(p-1)\dx &x_p=p\dx, &1\le p \le \nsx,\\
y_{q-1}=(q-1)\dy &y_q=q\dy, &1\le q \le \nsy,  \\
z_{r-1}=(r-1)\dz &z_r=r\dz,&1\le r \le \nbz.
\end{array}
\end{equation}
 The geometry is illustrated in Figure~\ref{figure1}, in which  the configuration of station at location $(i,j)$ relative to volume prism $pqr$  is shown.  
Here  the blocks in the $z$-direction define depth $z\ge0$ pointing down. Thus, without any padding around the domain, we assume there is one station located above each prism, so that there are $\nsx$ and $\nsy$ blocks in the $x$ and $y$-directions, respectively. 

We suppose that the entries in $G$ depend on the integral of kernel $h$ and correspond to the unit contribution from a given prism to a particular station. The ordering of the entries depends on the  organization of the volume domain into a vector of length, $n=\nbx\nby\nbz$. We will assume the multilayer model in which the volume is organized  by slices in depth and yields
\begin{equation}\label{Gslice}
G=[G^{(1)}, G^{(2)}, \dots, G^{(\nbz)}].
\end{equation}
Each $\Gr$ has size $\nsx\nsy \times  \nbx\nby =m\times \nbr$, where there are $\nbr$ prisms in slice $r$,  and maps from  the prisms in depth slice $r$, with depth coordinates $z_{r-1}$ and $z_{r}$, to the station measurements. Further, $\Gr$ decomposes as a block matrix with  block entries $\Gr_{jq}$, $1\le j \le \nsy$, $1\le q \le \nby$ each of size $\nsx\times\nbx$. Equivalently, this means that a given slice of the volume with $\nsy\nby$ blocks is mapped to a one dimensional vector using  \textbf{row-major} ordering; we sweep through the prisms in the slice for increasing  $x$ and fixed $y$ direction.  Entry $(\Gr)_{k\ell}$, $1\le k \le \nsx\nsy$,  $1\le \ell \le \nbx\nby$ represents the contribution from the prism at location $\ell=(q-1)\nbx+p$, for $1\le q \le \nby$ and $1\le p \le \nbx$, for depth slice $r$, to the station at $k=(j-1)\nsx+i$, $1\le j \le \nsy$, $1\le i \le \nsx$. We use $\tilde{h}(s_{ij})_{pqr}$ to denote the function that calculates the contribution to station $s_{ij}$ from prism $c_{pqr}$. Then
\begin{equation}\label{htilde}(\Gr)_{k\ell}=\tilde{h}(s_{ij})_{pqr}, \,\, k=(j-1)\nsx+i, \,\, \ell=(q-1)\nbx+p.
\end{equation}
Assuming that the discussion is applied for slice $r$,  we remove the dependence of $\tilde{h}$ on depth and use $\tilde{h}(s_{ij})_{pq}$ to indicate the contribution to station $s_{ij}=(a_{ij},b_{ij})$ due to block number $p$ in $x$ and $q$ in $y$. We note, further, that while the discussion is applied under the assumption of a uniform depth interval, $\dz$, the approach applies equally well for problems in which the multilayer coordinate grid has layers of different depths, see e.g. \cite{ZhangWong:15}.

\subsection{Spatially invariant kernels}\label{sec:invariantkernel}
Our discussion focuses on kernels that are spatially invariant in all dimensions: 
$$h(a,b,c,x,y,z)=h(x-a,y-b,z-c).$$ Then, again considering a single slice at depth $z$, 
the calculation of \eqref{htilde} depends on  the differences $(x-a)$ and $(y-b)$ for all station and prism coordinates, \eqref{stationsij} and \eqref{prismspqr}, respectively. Using the matrices
 \begin{align*}\left.\begin{array}{cc}(DX)_{ij,p}=(x_{p}-a_{ij})&0\le p \le \nbx  \\ (DY)_{ij,q}=(y_{q}-b_{ij})&0\le q \le \nby\end{array} \right\} 1\le i \le \nsx,\,\, 1\le j \le \nsy,\end{align*} 
the distances for block $pq$, $1\le p \le \nbx$ and $1\le q \le \nby$,  are obtained from distance matrices $(DX)_{p-1}$ and $(DX)_{p}$ in $x$ and likewise from $(DY)_{q-1}$ and $(DY)_q$ in $y$.  Now, these matrices are independent of the slice coordinate $r$, and, under the assumption that  the prisms are uniform in the $x-$ and $y-$ dimensions, 
 \begin{align*}
 (DX)_{p}=(DX)_{p-1}+\dx, \,\,1\le p \le \nbx, \,\,(DY)_{q}=(DY)_{q-1}+\dy, \,\,1\le q \le \nby. \end{align*}
Thus, all matrices  $(DX)_p$ and $(DY)_q$ can be obtained directly from $(DX)_0$ and $(DY)_0$ and they are independent of the slice, regardless of the locations of the stations relative to the prisms. 
When the stations are on a uniform grid so that $a_{ij}$ is independent of $j$ and $b_{ij}$ is independent of $i$, then the sizes of matrices $(DX)_0$ and $(DY)_0$ are reduced in the first dimension to  $\nsx$ and $\nsy$, respectively. We note that it is  practical, therefore,  to store $(DX)_0$ and $(DY)_0$ entirely, and  update an entire slice of the domain without recalculating $(DX)_0$ and $(DY)_0$ across slices,  regardless of the station locations. Still,  greater optimization is achieved when the stations are located also on a uniform grid. 

\subsubsection{Placement of the stations at the center of the cells}\label{sec:stationplacement}
Now, following \cite{BoCh:2001,ChenLiu:18}, suppose that the stations are placed at the centers of the cells such that $a_i=  (i-\frac12)\dx$, $1\le i \le \nsx$ and  $b_j=(j-\frac12)\dy$, $1\le j \le \nsy$. Then, the two coordinate systems for the stations and the volume domain, are uniformly staggered in the $x-y$ plane. Thus the distances between stations and coordinates are uniform, 
\begin{align*}
\begin{array}{lllll}
(DX)_{i,p}&=x_{p}-a_i &= (p-1) \dx -(i-\frac12) \dx &= (p-i-\frac12) \dx, & 1\le p \le \nbx+1 \\ 
(DY)_{j,q}&=y_{q}-b_j &=(q-1) \dy -(j-\frac12) \dy  &= (q-j-\frac12) \dy, & 1\le q \le \nby+1, 
\end{array}
\end{align*}
and for all pairs of indices $(i,p)$ and $(j,q)$, the possible paired distances are obtained from the vectors
\begin{align}\label{Distances}
\begin{array}{llll}
X_{\ell}&= (\ell-\nsx -\frac12)\dx, & 1\le \ell \le 2\nsx \\ 
Y_k&=(k-\nsy-\frac12) \dy,  &1\le  k \le 2\nsy. 
\end{array}
\end{align}
\subsubsection{Introducing padding around the domain}\label{sec:stationplacementpadding}
\begin{figure}
\begin{center}
\includegraphics[width=.75\textwidth]{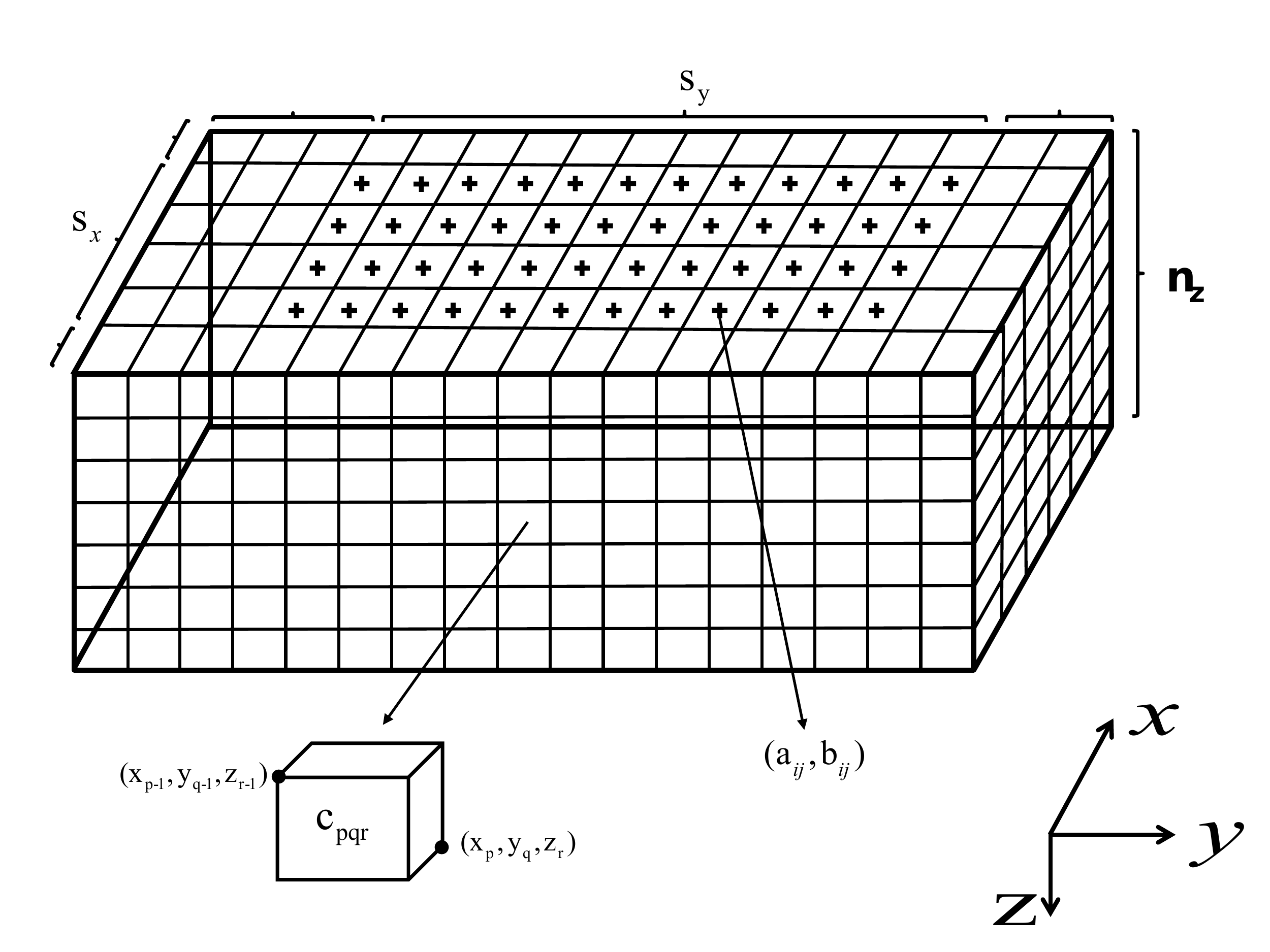}
\caption{The configuration of the volume domain with padding. \label{figure2}} \end{center}
\end{figure}

Suppose now that padding is introduced around the domain, with an extra $\padxl$ and $\padxr$  blocks in the $x$-direction,  so that the $x-$coordinates  extend from   $(-\padxl:(\nsx+\padxr))\dx$ for a total of $\nsx$ coordinate blocks within the domain but a total number of blocks $\nbx=(\nsx+\padxl+\padxr)$, where blocks $1$ to $\padxl$ are in the padded region to the left of the domain, and blocks $\nsx+\padxl+1$ to $\nbx$ are within the padded region to the right. Thus, the coordinates of block $p$ are adjusted to $(p-\padxl-1)\dx$ to $(p-\padxl)\dx$, consistent with \eqref{prismspqr} for $\padxl=0$. Likewise, the $y$ coordinates extend from $-\padyl:(\nsy+\padyr)$ and $\nby=(\nsy+\padyl+\padyr)$, see Figure~\ref{figure2}. Hence, 
 \eqref{Distances} is replaced   by 
\begin{align}\label{PaddedDistances}
\begin{array}{llll}
X_{\ell}&= (\ell-(\nsx+\padxl) -\frac12)\dx, & 1\le \ell \le 2\nsx+\padxl+\padxr =\nbx+\nsx\\ 
Y_k&=(k-(\nsy+\padyl)-\frac12) \dy,  &1\le  k \le 2\nsy+\padyl+\padyr =\nby+\nsy. 
\end{array}
\end{align}

If all possible paired distances are needed for calculating $\tilde{h}$, we only need to store vectors $X$ and $Y$, as given by \eqref{PaddedDistances}, which is negligible as compared to the entire storage of the $mn$ entries in the matrix $G$. Note that \eqref{PaddedDistances} explicitly assumes no stations in the padded regions. 

\subsection{Matrix Structure  for spatially invariant  kernels}\label{sec:kernels}
We now consider the structure of the matrices that arise for spatially invariant kernels,  and the associated computations that are required. We present the discussion first for domains without padding in Section~\ref{sec:symBTTBstructure}, and then the modifications that are required when domain padding is introduced, Section~\ref{sec:symBTTBstructurepadded}. The summary of the discussion is detailed in the presented Algorithms in  \ref{App:A}, and can be ignored if the intent is to only use the provided codes. 

\subsubsection{Symmetric kernel matrices with Toeplitz block structure without domain padding}\label{sec:symBTTBstructure}
We consider first the case in which the matrix $\Gr$ is  symBTTB.  Specifically, we suppose that the slice matrix $\Gr$ has a symmetric block structure and is defined by its first row $\Gr_{1q}=\Gr_q$, $1\le q \le \nby$. Then, with $\nbx=\nsx$ and $\nby=\nsy$, 
\begin{align}\label{symBTTBmatrix}
\Gr=\left[\begin{array}{ccccc} \Gr_1&\Gr_2&\Gr_3&\hdots &\Gr_{\nby} \\\Gr_2 & \Gr_1& \Gr_2 &\hdots &\Gr_{\nby-1} \\ \vdots & \ddots&\ddots&\ddots&\vdots \\ \vdots & \ddots&\ddots&\ddots&\vdots \\ \Gr_{\nsy} &\Gr_{\nsy-1}& \Gr_{\nsy-2}&\hdots&\Gr_1 \end{array}\right], \end{align}
where each $\Gr_q$ is  symmetric and defined by its first row, 
$$\Gr_q=\left[\begin{array}{ccccc}  g_{1q}&g_{2q}& g_{3q} &\hdots &g_{\nbx q} \\g_{2q} & g_{1q}& g_{2q} &\hdots &g_{(\nbx -1) q} \\ \vdots & \ddots&\ddots&\ddots&\vdots \\ \vdots & \ddots&\ddots&\ddots&\vdots \\ g_{\nsx q} &g_{(\nsx -1) q}& g_{(\nsx -2) q}&\hdots&g_{1q} \end{array}\right].$$
 MATLAB notation can be used to write these matrices compactly in terms of the defining first row (column). Specifically,  using MATLAB notation, 
\begin{align}\label{symmetrictoeplitzmatrix}\Gr_q=\mathrm{toeplitz}(\bfr_q), \,\,\bfr_q= (g_{1q}, g_{2q},  g_{3q}, \hdots ,g_{\nbx q} ),
\end{align} 
and,  with abuse of the same notation as applied to matrices, 
\begin{align}\label{symmetrictoeplitzblockmatrix}\Gr=\mathrm{toeplitz}(R), \,\,R= (\Gr_{1}, \Gr_{2},  \Gr_{3}, \hdots ,\Gr_{\nby}).
\end{align}

From \eqref{symmetrictoeplitzmatrix} and \eqref{symmetrictoeplitzblockmatrix} it is immediate, as discussed in \cite{BoCh:2001} and \cite{ChenLiu:18}, that the generation of $\Gr$ requires only the calculation of its first row. But the first row represents the contributions of all prisms to the first station. Thus, to find $\Gr$ requires only the calculation of  
\begin{align}\label{symtoeplitz}(\Gr)_{1\ell}&=\tilde{h}(s_{11})_{pq}, \,\, \ell=(q-1)\nbx+p, \,\, 1\le p \le \nbx, \,\, 1\le q \le \nby \quad \text{or}\\
(\Gr_1)_{1:\nbx\nby}&=(\begin{array}{c|c|c|c}\bfr_1&\bfr_2& \dots &\bf r_{\nby}
\end{array}) \quad \text{where} \nonumber  \\
\bfr_q&=(\tilde{h}(s_{11})_{1q}, \tilde{h}(s_{11})_{2q},  \tilde{h}(s_{11})_{3q}, \hdots ,\tilde{h}(s_{11})_{\nbx q} ), \,\, 1\le q \le \nby.\label{rqsymtoeplitz}
\end{align}
Equivalently, it is sufficient to calculate only the distances $(DX)_{1,p}=(p-3/2)\dx$ and $(DY)_{1,q}=(q-3/2)\dy$, for 
$1\le p \le \nbx+1$ and $1\le q \le \nby+1$, and \eqref{Distances} is replaced by 
\begin{align}\label{symXY}
\begin{array}{ll}
X_\ell=(\ell-\frac32) \dx, & 1\le  \ell \le (\nbx+1)\\ 
Y_k=(k-\frac32)  \dy, & 1\le  k \le (\nby+1).%
\end{array}
\end{align}

\subsubsection{Symmetric kernel matrices with Toeplitz block structure and domain padding}\label{sec:symBTTBstructurepadded}
We now consider what will happen when the domain is padded in the $x$ and $y$ directions, with no real stations within the padded region. To illustrate the impact of the padding on the generation of the matrix we take a one-dimensional example with $\nsx=4$, $\padxl=2$ and $\padxr=1$. Suppose first that there are artificial stations in the first two blocks and in the final block, namely for blocks $1$, $2$, and  $7$. Then the single square and symmetric Toeplitz that defines $\Gr$ is
\begin{align*}
\Gr_1=\left[\begin{array}{cc|cccc|c}
g_1&g_2&g_3&g_4&g_5&g_6&g_7\\
g_2&g_1&g_2&g_3&g_4&g_5&g_6\\ \hline
g_3&g_2&g_1&g_2&g_3&g_4&g_5\\
g_4&g_3&g_2&g_1&g_2&g_3&g_4\\
g_5&g_4&g_3&g_2&g_1&g_2&g_3\\
g_6&g_5&g_4&g_3&g_2&g_1&g_2\\ \hline
g_7&g_6&g_5&g_4&g_3&g_2&g_1
\end{array}\right] \quad 
\begin{array}{|l}\text{Station }1 (\text{Artificial}) \\ \text{Station }2 (\text{Artificial}) \\ \hline  \text{Station } 3=\padxl+1\\ \text{Station } 4 \\ \text{Station } 5 \\ \text{Station } 6=\padxl+\nsx \\ \hline \text{Station } 7 (\text{Artificial}) \end{array}.
\end{align*}
This depends on $$\bfr=(g_1, g_2, g_3, g_4, g_5, g_6, g_7)=(\tilde{h}(s_{1})_{1}, \tilde{h}(s_{1})_{2},  \tilde{h}(s_{1})_{3}, \hdots ,\tilde{h}(s_{1})_{7 } ).$$ But the contribution to the first real station due to all prisms is given by the third row, row $\padxl+1$ of $\Gr_1$, which is 
$$(\Gr_1)_{3}=\left(g_3,g_2,g_1,g_2,g_3,g_4,g_5 \right), $$
and the contributions for the  real stations  are determined, using symmetry,  by
\begin{align*}(\Gr_1)(\padxl+1:\padxl+\nsx, :) &=
\left[\begin{array}{ccccccc}
g_3&g_2&g_1&g_2&g_3&g_4&g_5\\
g_4&g_3&g_2&g_1&g_2&g_3&g_4\\
g_5&g_4&g_3&g_2&g_1&g_2&g_3\\
g_6&g_5&g_4&g_3&g_2&g_1&g_2\\
\end{array}\right] \\
&=\mathrm{toeplitz}(\bfc,\bfr). 
\end{align*}
Here 
\begin{align*}
 \bfc&=(g_3,g_4,g_5,g_6)\text{ and }  \bfr=( g_3,g_2,g_1,g_2,g_3,g_4,g_5). 
 \end{align*}
 More generally, for one dimension only, 
 \begin{align*}
 \bfc&=(g_{\padxl+1},g_{\padxl+2},\hdots ,g_{\padxl+\nsx})
= (\tilde{h}(s_{1})_{\padxl+1},\tilde{h}(s_{1})_{\padxl+2},\hdots ,\tilde{h}(s_{1})_{\padxl+\nsx})\text{ and }  \\
 \bfr& =( g_{\padxl+1},
 \hdots,g_2,g_1,g_2,\hdots,g_{\nbx-\padxl})\\
& =(\tilde{h}(s_{1})_{\padxl+1},\hdots,\tilde{h}(s_{1})_{2},\tilde{h}(s_{1})_{1}, \tilde{h}(s_{1})_{2},\hdots ,\tilde{h}(s_{1})_{\nbx-\padxl }).
\end{align*}
Extending to the two-dimensional case, and assuming that the first artificial station is in the $(1,1)$ block  of the padded domain,  then $\Gr_q$, for any $q$,  is also Toeplitz and is given by
\begin{align}\label{paddedsymtoeplitz} 
\Gr_q&=\mathrm{toeplitz}(\bfc_q,\bfr_q),\,\, 1\le q \le \nby, \\ \label{cqpaddedsymtoeplitz}
 \bfc_q&=(\tilde{h}(s_{11})_{(\padxl+1)q},\tilde{h}(s_{11})_{(\padxl+2)q},\hdots ,\tilde{h}(s_{11})_{(\padxl+\nsx)q})\text{ and }\\
 \bfr_q&=(\tilde{h}(s_{11})_{(\padxl+1)q},\hdots,\tilde{h}(s_{11})_{2q},\tilde{h}(s_{11})_{1q}, \tilde{h}(s_{11})_{2q},\hdots, 
 \tilde{h}(s_{11})_{(\nbx-\padxl)q} ).\label{rqpaddedsymtoeplitz}
\end{align}
This is consistent with \eqref{symtoeplitz} - \eqref{rqsymtoeplitz} for the unpadded case. But notice, also, that the maximum distance between station and coordinates  in the  $x-$coordinate is $\max(\nbx -\padxl, \nbx-\padxr )\dx=\max(\nsx+\padxr,\nsx+\padxl)\dx$. 

It remains to apply the same argument to the structure of the matrix $\Gr$, as to the structure of its individual components, to determine the structure of the symBTTB matrix when padding is applied. Then, consistent with \eqref{paddedsymtoeplitz}-\eqref{rqpaddedsymtoeplitz},  \eqref{symmetrictoeplitzblockmatrix} is replaced by 
\begin{align}\label{paddedsymmetrictoeplitzblockmatrix}
\Gr&=\mathrm{toeplitz}(C,R),  \\ \label{Cpaddedsymtoeplitz}
C&= (\Gr_{\padyl+1}, \hdots, \Gr_{\padyl+\nsy})\text{ and } \\
R&= (\Gr_{\padyl+1}, \hdots ,\Gr_2,\Gr_{1}, \Gr_{2},  \Gr_{3}, \hdots ,\Gr_{\nby-\padyl}).\label{Rpaddedsymtoeplitz}
\end{align} 
Moreover, since this matrix depends on the first row of the symmetric matrix, defined with respect to the artificial station at $s_{11}$, it is sufficient to still use \eqref{symXY} for the calculation of the relevant distances between the first station and all coordinate blocks. But, from \eqref{cqpaddedsymtoeplitz} - \eqref{rqpaddedsymtoeplitz}, and \eqref{Cpaddedsymtoeplitz}-\eqref{Rpaddedsymtoeplitz},  just as we do not calculate all entries $g_j$ in $\Gr_q$, we also do not calculate all the blocks $\Gr_q$, rather the blocks needed are for $q=1:\max(\nby-\padyl, \nby-\padyr)$. Thus, while we may use \eqref{symXY} to calculate the relevant distances, in practice some savings in memory and computation can be made, when padding is significant relative to $\nsx$ and $\nsy$,  by using 
\begin{align}\label{shortsymXY}
\begin{array}{ll}
X_\ell=(\ell-\frac32) \dx, &  1 \le \ell  \le \nsx+\max(\padxr, \padxl)+1 \\ 
Y_k=(k-\frac32)  \dy, &  1\le k \le \nsy+\max(\padyr, \padyl)+1  .%
\end{array}
\end{align} 

\subsubsection{Nonsymmetric kernel matrices with block structure}\label{sec:BTTBstructure} 
Consider now the non symmetric BTTB matrix  given by 
\begin{align}\label{nonsymBTTBmatrix}
\Gr=\left[\begin{array}{cccccc} \Gr_1&\Gr_2&\Gr_3&\hdots &\hdots&\Gr_{\nby} \\
\barGr_2 & \Gr_1& \Gr_2 &\hdots &\hdots&\Gr_{\nby-1} \\ 
\vdots & \ddots&\ddots&\ddots&\ddots&\vdots \\ 
\vdots & \ddots&\ddots&\ddots&\ddots&\vdots \\ 
\barGr_{\nsy} &\barGr_{\nsy-1}& \barGr_{\nsy-2}&\hdots&\hdots&\Gr_{1} \end{array}\right], \end{align}
where, without padding, we assume again that $\nby=\nsy$. 
This matrix depends on the first block row and column only, and, again using the abuse of the Toeplitz notation, is given by
\begin{align}\label{nonsymBTTBtoeplitz}
\Gr &= \mathrm{toeplitz}(C,R),\,\, C=(\Gr_1, \barGr_2,  \hdots ,\barGr_{\nsy}), \,\, R=(\Gr_1, \Gr_2,  \hdots ,\Gr_{\nsy}).
\end{align}
Here we use $\barGr_j$ to denote the contributions below the diagonal, and $\Gr_q$ for the contributions above the diagonal. None of the block matrices are symmetric and, therefore, to calculate $\Gr$ it is necessary to calculate  columns and rows that define the entries in $C$ and $R$. Calculating $R$ uses just the first row entries $\Gr_q$, but since each of these is not symmetric we need  also the first columns $\bfc_q$ of each block in $\Gr_q$. 

Using \eqref{htilde}, the $\Gr_q$ are given by \eqref{paddedsymtoeplitz} with
\begin{align}\label{toeplitzmatrixnonsymr}
\bfr_q&=(\tilde{h}(s_{11})_{1q}, \tilde{h}(s_{11})_{2q}, \hdots, \tilde{h}(s_{11})_{\nbx q}),\,\, 1\le q \le \nby \\
\bfc_q&=( \tilde{h}(s_{11})_{1q},\tilde{h}(s_{21})_{1q},\hdots, \tilde{h}(s_{\nsx 1})_{1q}), \,\, 1\le q \le \nby.\label{toeplitzmatrixnonsymc}
\end{align}
Effectively, rather than calculating  all entries in the first block row of $\Gr$, $(\nsx\nbx)\nby$ entries,  for each matrix of the block we calculate just its first row and column, for a total of $(\nsx+\nbx)\nby$ entries.

This leaves the calculation of the $\bar{G}^{(r)}_j$, $2\le j \le \nsy$, which by the Toeplitz structure only use the entries of the first block column of $\Gr$. They are given by
\begin{align}\barGr_j&=\mathrm{toeplitz}(\barc_j,\barr_j), \,\,  2\le j \le \nsy, 
\label{toeplitzbarmatrixnonsym}
\\ \nonumber 
\barc_j&=(\tilde{h}(s_{1j})_{11}, \tilde{h}(s_{2j})_{11}, \hdots, \tilde{h}(s_{\nsx  j})_{11}),\,\, 2\le j \le \nsy, \text{ and }\\ \label{toeplitzbarmatrixnonsymc}
\barr_j&=(\tilde{h}(s_{1j})_{11}, \tilde{h}(s_{1j})_{21},\hdots, \tilde{h}(s_{1j})_{\nbx 1}), \,\,2\le j \le \nsy.
\end{align}
It is also immediately clear that this requires not only all distances between the first station and all prism coordinates, as in \eqref{symXY}, but also for all stations and the first coordinate block (for the first column of $\Gr$) which  uses $X_{\ell}=(\ell+\frac12) \dx$,  $-\nsx \le \ell \le 0$, 
 and likewise for $Y_k$. Thus,  we require the full set of differences \eqref{Distances}.

 \subsubsection{Nonsymmetric kernel matrices with block structure and domain padding}\label{sec:BTTBstructurepadded} 
 As for the discussion of the impact of the domain padding on the symmetric kernel in Section~\ref{sec:symBTTBstructurepadded}, we first present an example  using one dimension, namely fixed $y$-coordinate. We again assume  $\nsx=4$, $\padxl=2$ and $\padxr=1$ and  suppose that there are artificial stations in the first two blocks and in the final block, namely for blocks $1$, $2$, and  $7$. Then, the single square but non-symmetric Toeplitz matrix that defines $\Gr$ is
\begin{align*}
\Gr_1=\left[\begin{array}{cc|cccc|c}
g_1&g_2&g_3&g_4&g_5&g_6&g_7\\
\gamma_2&g_1&g_2&g_3&g_4&g_5&g_6\\ \hline
\gamma_3&\gamma_2&g_1&g_2&g_3&g_4&g_5\\
\gamma_4&\gamma_3&\gamma_2&g_1&g_2&g_3&g_4\\
\gamma_5&\gamma_4&\gamma_3&\gamma_2&g_1&g_2&g_3\\
\gamma_6&\gamma_5&\gamma_4&\gamma_3&\gamma_2&g_1&g_2\\ \hline
\gamma_7&\gamma_6&\gamma_5&\gamma_4&\gamma_3&\gamma_2&g_1
\end{array}\right] \quad 
\begin{array}{|l}\text{Station } 1 (\text{Artificial})\\ \text{Station } 2 (\text{Artificial}) \\ \hline  \text{Station } 3 =\padxl+1\\ \text{Station } 4 \\ \text{Station } 5\\ \text{Station } 6=\padxl+\nsx \\ \hline \text{Station } 7 (\text{Artificial})\end{array}
\end{align*}
This depends on 
\begin{align*}
\bfr&=(g_1, g_2, g_3, g_4, g_5, g_6, g_7)=(\tilde{h}(s_{1})_{1}, \tilde{h}(s_{1})_{2},   \hdots ,\tilde{h}(s_{1})_{6},\tilde{h}(s_{1})_{7 } )\\
\bfc&=(g_1, \gamma_2, \gamma_3, \gamma_4, \gamma_5, \gamma_6, \gamma_7)=(\tilde{h}(s_{1})_{1}, \tilde{h}(s_{2})_{1}, \hdots,
\tilde{h}(s_{6})_{1},  \tilde{h}(s_{7})_{1 } ).
\end{align*}
But again the  required rows of $\Gr_1$ are those that correspond to the actual real stations
\begin{align*}(\Gr_1)(\padxl+1:\padxl+\nsx, :) &=
\left[\begin{array}{ccccccc}
\gamma_3&\gamma_2&g_1&g_2&g_3&g_4&g_5\\
\gamma_4&\gamma_3&\gamma_2&g_1&g_2&g_3&g_4\\
\gamma_5&\gamma_4&\gamma_3&\gamma_2&g_1&g_2&g_3\\
\gamma_6&\gamma_5&\gamma_4&\gamma_3&\gamma_2&g_1&g_2\\
\end{array}\right] \\
&=\mathrm{toeplitz}(\bfc,\bfr) \text{ where } \\
 \bfc&=(\gamma_3,\gamma_4,\gamma_5,\gamma_6)\text{ and }  \bfr=( \gamma_3,\gamma_2,g_1,g_2,g_3,g_4,g_5). 
 \end{align*}
 More generally, 
 \begin{align*}
 \bfc&=(\gamma_{\padxl+1},\gamma_{\padxl+2},\hdots ,\gamma_{\padxl+\nsx})= (\tilde{h}(s_{\padxl+1})_{1},\tilde{h}(s_{\padxl+2})_{1},\hdots ,\tilde{h}(s_{\padxl+\nsx})_{1})\text{ and }  \\
 \bfr& =( \gamma_{\padxl+1},\gamma_{\padxl},\hdots, \gamma_2,g_1,g_2,\hdots,g_{\nbx-\padxl})\\
 &=(\tilde{h}(s_{\padxl+1})_{1},\tilde{h}(s_{\padxl})_{1}, \hdots,\tilde{h}(s_{2})_{1},\tilde{h}(s_{1})_{1}, \tilde{h}(s_{1})_{2},\hdots ,\tilde{h}(s_{1})_{\nbx-\padxl }).
\end{align*}
Extending to the two-dimensional case, with the same assumptions as in Section~\ref{sec:symBTTBstructurepadded}, $\Gr_q$ is obtained as 
\begin{align}\label{paddedtoeplitz} 
\Gr_q&=\mathrm{toeplitz}(\bfc_q,\bfr_q),\,\, 1\le q \le \nby, \\ 
 \bfc_q&=(\tilde{h}(s_{(\padxl+1)1})_{1q},\tilde{h}(s_{(\padxl+2)1})_{1q},\hdots ,\tilde{h}(s_{(\padxl+\nsx)1})_{1q})\text{ and } \nonumber\\
 \bfr_q&=(\tilde{h}(s_{(\padxl+1)1})_{1q},\tilde{h}(s_{\padxl 1})_{1q},\hdots,\tilde{h}(s_{21})_{1q} 
 \tilde{h}(s_{11})_{1q},\hdots, \tilde{h}(s_{11})_{(\nbx-\padxl)q} ).\nonumber
\end{align}
This is consistent with \eqref{toeplitzmatrixnonsymr} - \eqref{toeplitzmatrixnonsymc} for the unpadded case. 

Turning to the column block entries, first observe that $\barGr_1=\Gr_1$, and so  we examine
$\barGr_j$ which represents stations $1$ to $\nbx$ (both real and artificial) in the $x$-direction for a fixed $j$ coordinate in the $y$-direction. Then, with the same example for choices of $\nsx$, $\padxl$ and $\padxr$, 
\begin{align*}
\barGr_j=\left[\begin{array}{cc|cccc|c}
\barg_1&\barg_2&\barg_3&\barg_4&\barg_5&\barg_6&\barg_7\\
\bargamma_2&\barg_1&\barg_2&\barg_3&\barg_4&\barg_5&\barg_6\\ \hline
\bargamma_3&\bargamma_2&\barg_1&\barg_2&\barg_3&g_4&\barg_5\\
\bargamma_4&\bargamma_3&\bargamma_2&\barg_1&\barg_2&\barg_3&\barg_4\\
\bargamma_5&\bargamma_4&\bargamma_3&\bargamma_2&\barg_1&\barg_2&\barg_3\\
\bargamma_6&\bargamma_5&\bargamma_4&\bargamma_3&\bargamma_2&\barg_1&\barg_2\\ \hline
\bargamma_7&\bargamma_6&\bargamma_5&\bargamma_4&\bargamma_3&\bargamma_2&g_1
\end{array}\right] \quad 
\begin{array}{|l}\text{Station } 1 (\text{Artificial})\\ \text{Station } 2 (\text{Artificial}) \\ \hline  \text{Station } 3 =\padxl+1\\ \text{Station } 4 \\ \text{Station } 5\\ \text{Station } 6=\padxl+\nsx \\ \hline \text{Station } 7 (\text{Artificial})\end{array}
\end{align*} 
This depends on 
\begin{align*}
\barr&=(\barg_1, \barg_2, \barg_3, \barg_4, \barg_5, \barg_6, \barg_7)=(\tilde{h}(s_{1j})_{11}, \tilde{h}(s_{1j})_{21}, \hdots,
\tilde{h}(s_{1j})_{61 },  \tilde{h}(s_{1j})_{71 } )\\
\barc&=(\barg_1, \bargamma_2, \bargamma_3, \bargamma_4, \bargamma_5, \bargamma_6, \bargamma_7)=(\tilde{h}(s_{1j})_{11}, \tilde{h}(s_{2j})_{11},  \hdots, 
 \tilde{h}(s_{6j})_{11},  \tilde{h}(s_{7j})_{11 } ).
\end{align*}
But again, since stations $1$ to $2$  and $7$ are artificial, we only need 
\begin{align*}(\barGr_j)(\padxl+1:\padxl+\nsx, :) &=
\left[\begin{array}{ccccccc}
\bargamma_3&\bargamma_2&\barg_1&\barg_2&\barg_3&\barg_4&\barg_5\\
\bargamma_4&\bargamma_3&\bargamma_2&\barg_1&\barg_2&\barg_3&\barg_4\\
\bargamma_5&\bargamma_4&\bargamma_3&\bargamma_2&\barg_1&\barg_2&\barg_3\\
\bargamma_6&\bargamma_5&\bargamma_4&\bargamma_3&\bargamma_2&\barg_1&\barg_2\\
\end{array}\right] \\
&=\mathrm{toeplitz}(\bfc,\bfr) \text{ where } \\
 \bfc&=(\bargamma_3,\bargamma_4,\bargamma_5,\bargamma_6)\text{ and }  \bfr=( \bargamma_3,\bargamma_2,\barg_1,\barg_2,\barg_3,\barg_4,\barg_5). 
 \end{align*}
Thus, in two dimensions, the first column block entries are $\barGr_j$, $1\le j \le \nby$, with 
\begin{align}\label{barpaddedtoeplitz}
\barGr_j&=\mathrm{toeplitz}(\barc_j,\barr_j), \,\, 1\le j \le \nsy+\padyl+\padyr, \nonumber\\ 
\barc_j&= (\tilde{h}(s_{(\padxl+1)j})_{11},\tilde{h}(s_{(\padxl+2)j})_{11},\hdots, \tilde{h}(s_{(\padxl+\nsx)j})_{11}) \text{ and }\\ \nonumber
\barr_j& = (\tilde{h}(s_{(\padxl+1)j})_{11}, \tilde{h}(s_{\padxl j})_{11},\hdots, \tilde{h}(s_{2j})_{11}, \tilde{h}(s_{1j})_{11},\hdots, \tilde{h}(s_{1j})_{(\nbx-\padxl)1}).
\end{align}
But now \eqref{nonsymBTTBtoeplitz} is replaced by the block Toeplitz matrix
\begin{align}\label{nonsymBTTBtoeplitzpadded}
\Gr&=\mathrm{toeplitz}(C,R),  \\ \nonumber 
C&= (\barGr_{\padyl+1}, \hdots, \barGr_{\padyl+\nsy}) \text{ and }\\ \nonumber
R&= (\barGr_{\padyl+1}, \barGr_{\padyl}, \hdots, \barGr_{2}, \Gr_1, \hdots, \Gr_{\nby-\padyl}).
\end{align}
Here each block matrix is the subset of rows corresponding to the real stations, as noted in \eqref{paddedtoeplitz} and \eqref{barpaddedtoeplitz}. Moreover, we conclude that \eqref{paddedtoeplitz} is applied only for $1\le q \le \nby-\padyl=\nsy+\padyr$ and \eqref{barpaddedtoeplitz} for $1 \le j \le \padyl+\nsy$, reducing the dimension of the required $Y$ in the $y-$direction. Likewise, $X$ is reduced because of the padding impacting the required entries for generating both $\Gr_q$ and $\barGr_j$. Thus while the required vectors are given by \eqref{Distances}, with $\nbx$ replacing $\nsx$ and $\nby$ replacing $\nsy$, their lengths can be reduced as  for the symmetric case, \eqref{shortsymXY}, by using
\begin{align}\label{shortunsymXY}
\begin{array}{ll}
X_\ell=(\ell-(\nsx+\max(\padxr, \padxl))-\frac12) \dx, & 1\le \ell \le 2(\nsx+\max(\padxr, \padxl)) \\ 
Y_k=(k-(\nsy+\max(\padyr, \padyl))-\frac12)  \dy, & 1\le k\le 2(\nsy+\max(\padyr, \padyl)) .%
\end{array}
\end{align} 
This effectively assumes the calculation of $\barGr_1$ as well as $\Gr_1$, whereas only one is calculated in practice, since $\barGr_1=\Gr_1$. 
 
\begin{figure}
 \begin{center}
 \includegraphics[width=.5\textwidth]{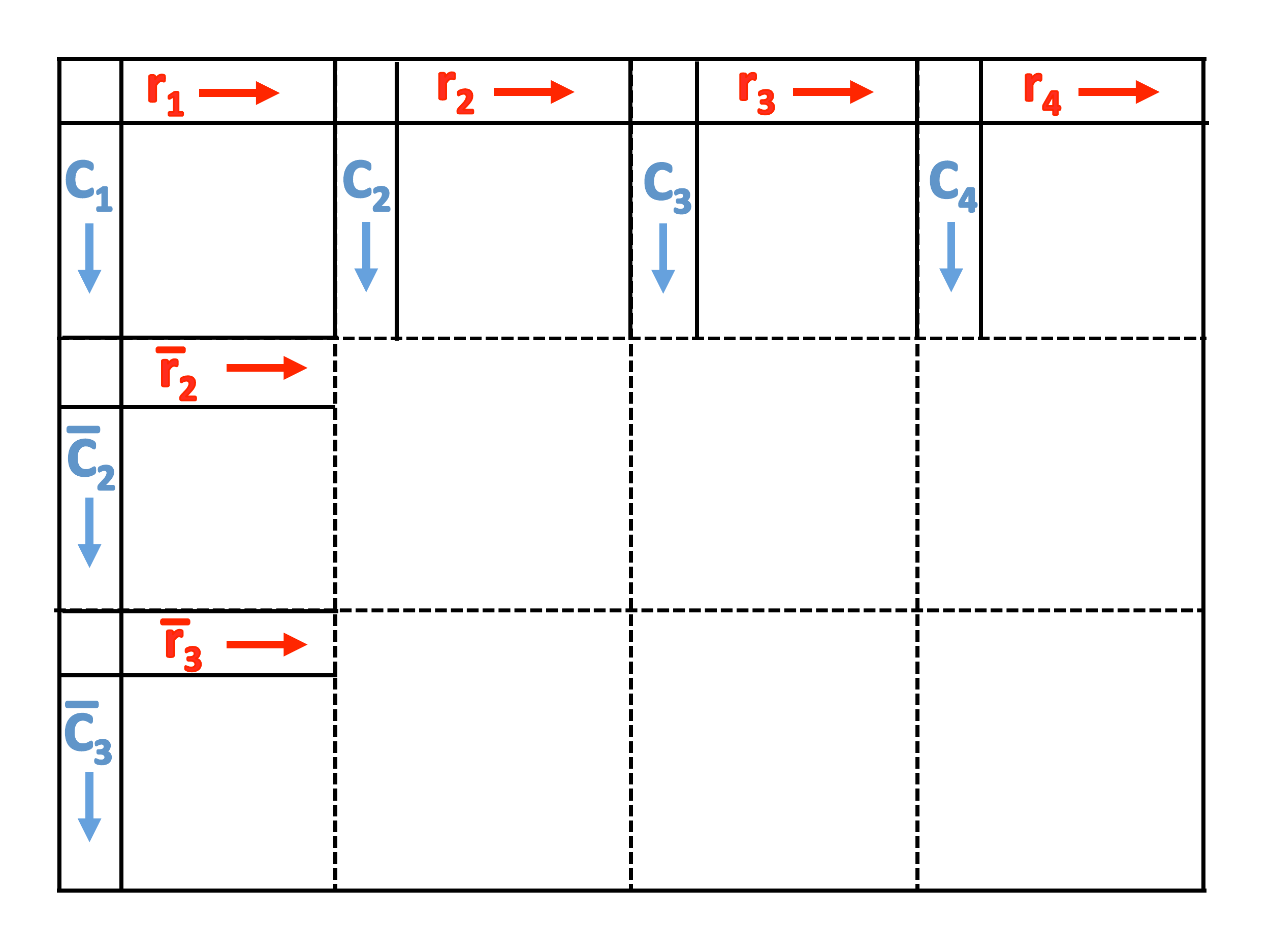}
 \caption{Required entries from a matrix $\Gr$ in order to calculate its complete block Toeplitz form. This shows considerable savings can be accrued in calculating the entries of $\Gr$ when using the structure. \label{figure3}}
 \end{center}  \end{figure}

The plot in  Figure~\ref{figure3} illustrates the unique entries from $\Gr$ that define its block Toeplitz structure. 

\section{Circulant Operators and the 2D FFT}\label{sec:circulant}
\begin{defn}[Circulant]\label{defn:circ}
The Toeplitz matrix in which the defining vectors $\bfc$ and $\bfr$, each of length $2\nsx-1$, have entries that are related by $r_{i} = c_{(2\nsx+1-i) }$ for $2\le i \le 2\nsx-1$, is \textbf{circulant}.  
\end{defn}
Notice that a circulant matrix is defined solely by its first column or first row.  Here we will use the first column. Moreover, any Toeplitz matrix can be embedded in its circulant extension, as illustrated for the simple example with $\nsx = \nbx = 3$ 
\begin{equation*}
\left[\begin{array}{ccc|cc}
g_1 & g_2 & g_3 & \gamma_3 & \gamma_2\\
\gamma_2 & g_1 & g_2 & g_3 & \gamma_3\\
\gamma_3 & \gamma_2 & g_1 & g_2 & g_3\\ \hline
g_3 & \gamma_3 & \gamma_2 & g_1 & g_2\\
g_2 & g_3 & \gamma_3 & \gamma_2 & g_1
\end{array}\right].
\end{equation*}
In the same way, a BTTB matrix can be embedded in a block circulant matrix. Thus, the matrices \eqref{symBTTBmatrix} and \eqref{nonsymBTTBmatrix} can be embedded in  block circulant matrices, in which also  each Toeplitz block $\Gr_q$ and $\barGr_q$ is embedded in a  $(2\nsx-1) \times (2\nsx-1)$ circulant matrix. This yields a matrix that is Block Circulant with Circulant Blocks (BCCB). It is the structure of a BCCB matrix that facilitates the use of the 2DFFT to efficiently evaluate forward matrix multiplication with a BTTB matrix. Specifically,  given a BTTB matrix, it has been shown in \cite{Vogel:2002} that matrix-vector multiplication can be applied at reduced computational cost  by using a BCCB extension combined with the FFT for implementing a  discrete convolution.  The required components that provide the FFT approach are now discussed.

\begin{defn}[Exchange matrix]\label{defn:exchange}
The \textbf{exchange} matrix is the $m\times m$  matrix $J_m$ which is everywhere $0$ except for  $1$'s on the principal counter diagonal.  \end{defn} 
Given arbitrary vector $\bfx$ of length $m$, with entries $x_i$, $1\le i \le m$,  then $J_m\bfx=\bfy$ where $\bfy_i=\bfx_{m-i+1}$, namely it is the  the vector with the order of the entries reversed. Equivalently, for matrix $A$ with rows $\bfa_i$, $1\le i \le m$, then $J_mA = B$ where $B$ is the matrix with rows in reverse order, $\bfb_i = \bfa_{m-i+1}$. Further, multiplying on the right reorders the columns in reverse order. Specifically, $J_m^T=J_m$, and thus $\bfy^T=(J_m\bfx)^T=\bfx^TJ_m^T=\bfx^T J_m$ and the column entries of $\bfy$ are in reverse order as compared to $\bfx$. In the same way, $A J_m$ gives the matrix with the columns in reverse order.  In MATLAB the exchange matrix is implemented using the functions $\flipud$ and $\fliplr$, for ``up-down" and  ``left-right", for multiplication with $J_m$ on the left and right, respectively.

The exchange matrix yields a compact notation for the entries that define the circulant extension of a Toeplitx matrix.  For matrix $\Gr_q$ which depends on  $\bfr_q$, as given in \eqref{symmetrictoeplitzmatrix}, then  the defining first row for the circulant extension for each symmetric $\Gr_q$ is  given by
\begin{equation}\label{Eqn:symCext}
\brext_q=\left(\begin{array}{c} \bfr_q\\ J_{\nsx-1}\bfr_q(2:\nsx) \end{array}\right).
\end{equation}
For the non-symmetric case for $\Gr_q$, as given in \eqref{toeplitzmatrixnonsymr}-\eqref{toeplitzmatrixnonsymc},  the circulant extension uses
\begin{equation}\label{Eqn:Cext}
\bcext_q=\left(\begin{array}{c} \bfc_q\\ J_{\nsx-1}\bfr_q(2:\nsx) \end{array}\right) \text{ and }
\brext_q=\left(\begin{array}{c} \bfr_q\\ J_{\nsx-1}\bfc_q(2:\nsx) \end{array}\right).
\end{equation}
An equivalent  expression applies for the circulant extension for each $\barGr_j$ as defined in \eqref{toeplitzbarmatrixnonsym}-\eqref{toeplitzbarmatrixnonsymc} using the extension for $\barc_j$ and $\barr_j$. While \eqref{Eqn:symCext} and \eqref{Eqn:Cext} can be used as the defining vectors to explicitly generate the extensions $(\Gr_j)^{\mathrm{circ}}$ and $(\barGr_j)^{\mathrm{circ}}$ as Toeplitz matrices, again using $\texttt{toeplitz}(\bcext_j, \brext_j)$, we note that the intent is to define the vectors that define the extensions but not to generate the extensions. Moreover, $\brext_j$ is as noted defined explicitly from $\bcext_q$ and we focus entirely on the columns $\bfc^{\mathrm{ext}}_j$.  
We also note that this definition for generating the extension differs from that used in \cite{Li2018,Vogel:2002,ZhangWong:15}; the extra $0$ is omitted for convenience.  We also directly define the circulant extension instead of performing a series of transformations.

We now turn to the defining set of vectors  needed for the circulant extension of \eqref{nonsymBTTBmatrix}, which depends on its first block row and column as given in \eqref{nonsymBTTBtoeplitz}. Using the same analogy as with the block Toeplitz matrix and using Definition~\ref{defn:circ}, the extension for $\Gr$ requires the extensions of $C$ and $R$ in \eqref{nonsymBTTBtoeplitz}, thus for entries $\Gr_j$ and $\barGr_j$ for $1\le j \le \nsx$. Moreover, the circulant extension, as in the one-dimensional case will depend entirely either on the extension of $C$ or $R$, denoted by $C^{\mathrm{ext}}$ and $R^{\mathrm{ext}}$, but again we do not form $\mathtt{toeplitz}(C^{\mathrm{ext}}, R^{\mathrm{ext}})$. Again we assume the use of the extension for $C$ only, and note that   $C^{\mathrm{ext}}$ is completely defined by $\Tcirc$, dropping the dependence on slice $r$. Then using \eqref{Eqn:Cext} applied for the matrix form \begin{equation}\label{Eqn:Tcirc}
T^{\mathrm{circ}} = \left(\begin{array}{cccccccc}\barcext_1 & \cdots & \barcext_{\nsy} & \bcext_{\nsy} &  \cdots & \bcext_2\\
\end{array}\right),
\end{equation}
 and is of size $(2\nsx-1) \times (2\nsy-1)$

Now, using $\bfu=\bvec(U)$ to denote the vectorization of matrix $U$,  the two-dimensional convolution product $\Gr \bfu$, can be computed using  $T^{\mathrm{circ}}$ which defines the circular extension of $\Gr$, \cite{Vogel:2002}. Specifically, suppose that $(\Gr)^{\mathrm{circ}}$ is defined by $\Tcirc$, and let $\bfw = \bvec(W)$. Then the reshaped convolution product $\barray((\Gr)^{\mathrm{circ}} \bfw)$, where $\barray()$ is the inverse of $\bvec()$,  can be computed by
\begin{equation}\label{GrviaT}
\barray((\Gr)^{\mathrm{circ}} \bfw)=\Tcirc \star W = \ifftt(\fftt(\Tcirc)\dotstar \fftt(W)).
\end{equation} 
Here $\star$ denotes convolution, $\fftt$ denotes the two-dimensional FFT, and $\ifftt$ denotes the inverse two-dimensional FFT, and we introduce $\hat{T}^{\mathrm{circ}}=\fftt(\Tcirc)$ and $\hat{W}=\fftt(W)$. But now to obtain $\Gr \bfu$ from this product we notice that  $\Gr$ is in the upper left block of $(\Gr)^{\mathrm{circ}}$.  Thus, we 
define $W$ of size ${(2\nsx-1) \times (2\nsy-1)}$ by
\begin{equation}\label{Eqn:W}
W = \left[\begin{array}{cc} U & 0_{\nsx(\nsy-1)}\\ 0_{(\nsx-1)\nsy} & 0_{(\nsx-1)(\nsy-1)}\end{array}\right],
\end{equation}
using $0_{mn}$ to denote a matrix of zeros of size $m \times n$, and 
with $U$ of size ${\nsx \times \nsy}$. Hence, 
$\barray(\Gr\bfu)$ is  the upper left $\nsx \times \nsy$ block of $\barray((\Gr)^{\mathrm{circ}}\bfw)$ in \eqref{GrviaT}.  Moreover, $\barray((\Gr)^{\mathrm{circ}})$ does not need to be formed explicitly for this product.  Instead we directly calculate the elements of $\Tcirc$ using \eqref{Eqn:Tcirc}.

It is immediate that a set of equivalent steps can be used to calculate $(\Gr)^T\bfv$, where $\bfv=\bvec(V)$ for matrix $V$, since $(\Gr)^T$ is also BTTB.  In addition, the defining first column for $((\Gr)^{\mathrm{circ}})^T$ is the first row of $(\Gr)^{\mathrm{circ}}$.  Extending the relationship of the first column and row of a circulant matrix to a BCCB matrix, the matrix defining the BCCB transpose is $\Tcirc$ with columns $2$ through $2\nsy -1$ swapped left to right, and rows $2$ through $2\nsx -1$ swapped top to bottom.  This can be achieved using the following two steps. First obtain $\Tcirc$ and let 
\begin{equation}\label{Eqn:Ttran0}
T = \left(\begin{array}{cc}\Tcirc(1:2\nsx-1,1) & \Tcirc( 1:2\nsx-1,2: 2\nsy-1)J_{2\nsy-2}\end{array}\right).
\end{equation}
Then the required matrix $\tilde{T}^{\mathrm{circ}} $, replacing $\Tcirc$ in \eqref{GrviaT} for the transpose operation, is given by
\begin{equation}\label{Eqn:Ttran}
\tilde{T}^{\mathrm{circ}}= \left(\begin{array}{c} T(1,1:2\nsy-1)\\ J_{2\nsx-2}T(2:2\nsx-1,1:2\nsy-1)\end{array}\right).
\end{equation}

\subsection{Convolution with Domain Padding}
Now suppose that padding is introduced around the domain and, consistent with \eqref{paddedsymtoeplitz} and \eqref{barpaddedtoeplitz}, we assume that the indices for $\bfr_p$ and $\bfc_q$ are from the first row and column of the padded domain.   Then, for the case of the symBTTB matrix, \eqref{Eqn:symCext} is replaced by

\begin{equation}\label{Eqn:symCextpadded}
\brext_q=\left(\begin{array}{c} \bfr_q\\ J_{\nsx-1}\bfc_q(2:\nsx) \end{array}\right), \quad 
\bcext_q=\left(\begin{array}{c} \bfc_q\\ J_{\nbx-1}\bfr_q(2:\nbx) \end{array}\right),
\end{equation}
where here we use the definitions  \eqref{cqpaddedsymtoeplitz} and \eqref{rqpaddedsymtoeplitz} for $\bfr_q$ and $\bfc_q$. But now since $\bfr_q$ is defined by $\bfc_q$ for the symmetric case we can use just $\bcext_q$. Each vector is of length $\nsx$ for $\bfc_q$ and  $\nbx-1$ for $J_{\nbx-1}(\bfr_q(2:\nbx))$. Then using \eqref{Cpaddedsymtoeplitz}, the BCCB extension is defined by the replacement of \eqref{Eqn:Tcirc} by the matrix
\begin{equation}
 \label{Eqn:Tcircpadded}
T^{\mathrm{circ}} = \left(\begin{array}{cccccc|ccc}\bcext_{1+\padyl} & \cdots & \bcext_{\nsy+\padyl} & 
 \bcext_{\nsy+\padyr}  & \cdots & \bcext_{2} &\bcext_{1} & \cdots & \bcext_{\padyl}\\
\end{array}\right),
\end{equation}
of size $(\nsx+\nbx-1)  \times (\nsy+\nby-1)   $. 
 For the nonsymmetric case \eqref{Eqn:Tcircpadded} is replaced by
\begin{equation}
 \label{Eqn:Tcircpaddedgeneral}
 T^{\mathrm{circ}} = \left(\begin{array}{cccccc|ccc}\barcext_{1+\padyl} &  \cdots & \barcext_{\nsy+\padyl} & 
 \bcext_{\nsy+\padyr} &  \cdots & \bcext_{2} &\barcext_{1} & \cdots & \barcext_{\padyl}\\
\end{array}\right),
\end{equation}
where $\barcext_j$ is obtained as in \eqref{Eqn:symCextpadded} but using $\barc_j$ and $\barr_j$ from \eqref{barpaddedtoeplitz}.
We note that in any case in which $\padyl=0$ then the end block is removed. Moreover, due to $\Gr \bfu$ of size $\nsx\nsy$, when $\bfu$ is of size $\nbx\nby$, the definition of $W$ in 
\eqref{Eqn:W} is replaced by
\begin{align}\label{Eqn:Wpadded}
W = \left[\begin{array}{cc} U & 0_{\nbx(\nsy-1)}\\ 0_{(\nsx-1)\nby} & 0_{(\nsx-1)(\nsy-1)}\end{array}\right],
\end{align}
with $U\in\mathcal{R}^{\nbx \times \nby}$.
 
The transpose operation can still be applied to obtain $\tilde{T}^{\mathrm{circ}}$ from  \eqref{Eqn:Ttran0}-\eqref{Eqn:Ttran}.  But notice that for  $\bfv$ of size ${\nsx\nsy}$,  $(\Gr)^T\bfv$ is of size ${\nbx\nby}$, and in this case 
 \eqref{Eqn:W} is replaced by
\begin{align}\label{Eqn:WTpadded}\tilde{W} = \left[\begin{array}{cc} V & 0_{\nsx(\nby-1)}\\ 0_{(\nbx-1)\nsy} & 0_{(\nbx-1)(\nby-1)}\end{array}\right],
\end{align}
with $V$ of size $\nsx \times \nsy$.

We illustrate in Figure~\ref{figure4} the matrix $\Tcirc$ that is generated using the row and column entries from the complete BTTB $\Gr$ as shown in Figure~\ref{figure3}.
\begin{figure}[ht!]\begin{center}
\includegraphics[width=0.45\textwidth]{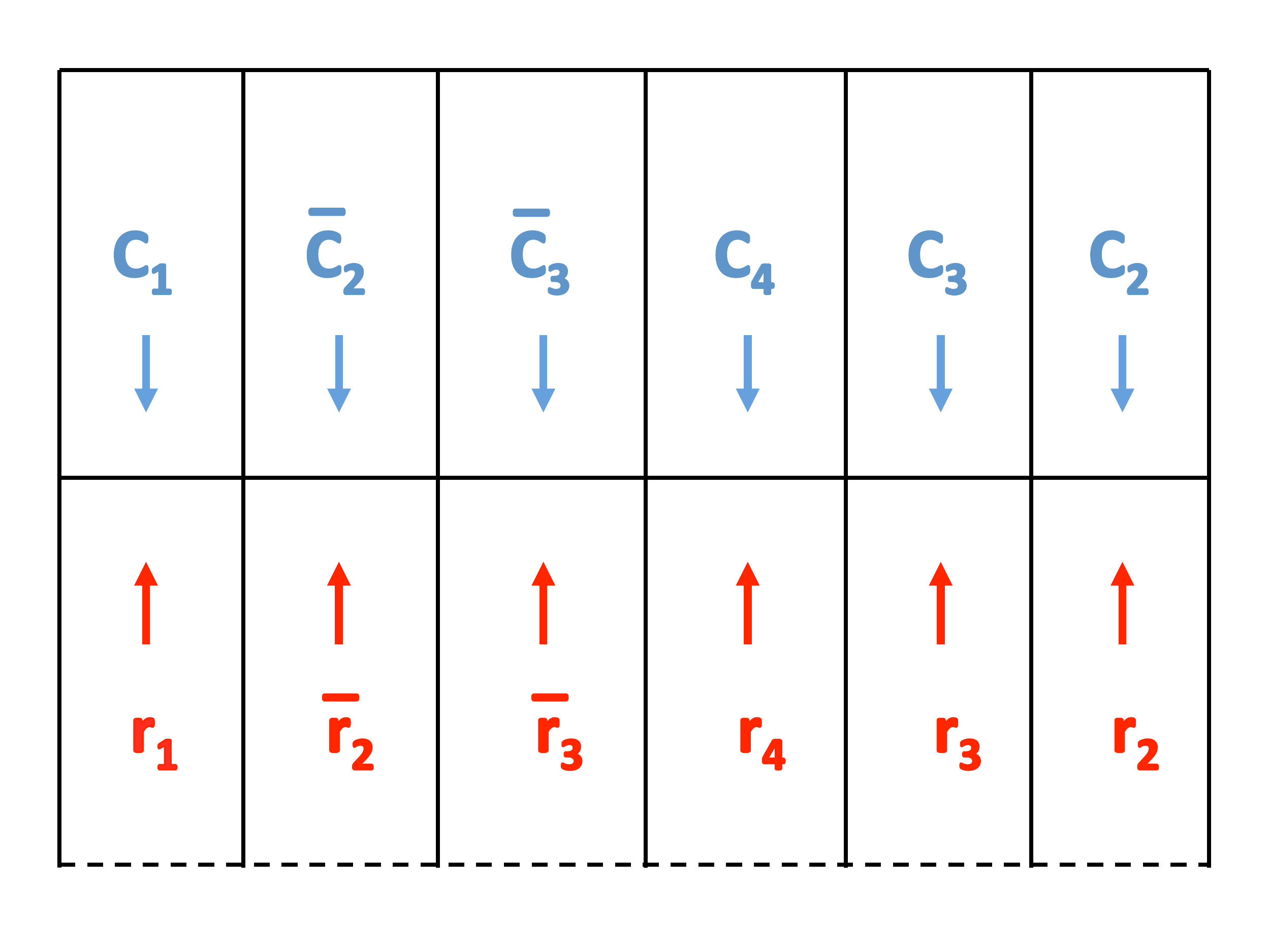}
\caption{The configuration of $\Tcirc$, where the arrow denotes the direction of the vector in ascending order. The   dotted line indicates that the first elements of each $\bm{r}_q$ and $\bar{\bm{r}}_q$ are omitted in the construction of $\Tcirc$. \label{figure4}}
\end{center}\end{figure}

\section{Optimizing the Calculations for Specific kernels}\label{sec:specifickernels}
\cite{ChenLiu:18} demonstrated that considerable savings are realized in the generation of the entries for the matrices $\Tcirc$ through optimized calculations of the entries necessary for finding $\Tcirc$ for forward modeling of the gravity problem. Here we focus on both improving that optimization and with the generation of an optimized and stable calculation for the magnetic kernel following the derivation of \cite{RaoBabu:91}.
\subsection{Gravity Kernel Calculation}\label{sec:gravitykernel}
The gravity kernel generates a symBTTB matrix for each slice in depth ($z$-direction).
According  to \cite{BoCh:2001}, and as used in \cite{ChenLiu:18}, the contribution of the kernel from the prism at point $(p,q)$ on the volume grid, (where $x$ points East and $y$ points North), to the station at location $(a,b,c)$, and here we assume  $c=0=z_0$,   is given by
\begin{align*}\tilde{h}(a,b,c)_{pq}&=\gamma \sum_{i=1}^2\sum_{j=1}^2\sum_{k=1}^2 (-1)^{i} (-1)^{j} (-1)^{k} \\
& \left(\txtzet_k \arctan \frac{\txtchi_i\txtups_j}{\txtzet_k \txtrho^k_{ij}}-\txtchi_i \ln \left(\txtrho^k_{ij}+\txtups_j\right)-\txtups_j \ln \left(\txtrho^k_{ij}+\txtchi_i\right)\right).\end{align*}
Here $\gamma$ is the gravitational constant and 
\begin{align}\label{notationgravity}
\begin{array}{ll}
\txtchi_1=x_{p-1} -a, &\txtchi_2=x_{p} -a \\
\txtups_1=y_{q-1} -b& \txtups_2=y_{q} -b\\
\txtzet_1=z_{r-1}-c,& \txtzet_2=z_{r}-c\\
R_{ij}^2 = \txtchi_i^2+\txtups^2_j &\txtrho^k_{ij}= \sqrt{R_{ij}^2+\txtzet_k^2}.
\end{array}
\end{align}
Specifically,  we need 
\begin{align*}&\sum_i \sum_j (-1)^{i+j+1}\left(\left(\txtzet_1 \arctan{\frac{(\txtchi\txtups)_{ij}}{\txtzet_1(R_1)_{ij}}} -\txtzet_2\arctan{\frac{(\txtchi\txtups)_{ij}}{\txtzet_2(R_2)_{ij}}}\right)\right.-\\ &\left.\txtchi_i \left(\ln((R_1)_{ij}+\txtups_j)-\ln((R_2)_{ij}+\txtups_j)\right) -\txtups_j \left(\ln((R_1)_{ij}+\txtchi_i)-\ln((R_2)_{ij}+\txtchi_i)\right)\right).\end{align*}
Here
$$R_1= \sqrt{R^{\wedge }2+\txtzet_1^2}, \quad R_2 =\sqrt{R^{\wedge }2+\txtzet_2^2},$$ and operations involve elementwise powers and multiplications. 
Using the notation in \cite{ChenLiu:18}, we write the summand, ignoring $\sum_j (-1)^{i+j+1}$ as 
\begin{align*}
\left(\left(\txtzet_1 (CM5)_{ij} -\txtzet_2(CM6)_{ij}\right) -\txtchi_i \left((CM3)_{ij}-(CM4)_{ij}\right) -\txtups_j \left((CM1)_{ij}-(CM2)_{ij}\right)\right).
\end{align*}
Now  notice that  $(CM3)_{ij}-(CM4)_{ij}$ is a logarithmic difference (and also for $(CM1)_{ij}-(CM2)_{ij}$). 
Thus, the differences can be replaced by 
\begin{equation}\label{simplifiedCMX} CMX=\ln{\frac{\txtchi+R_1}{\txtchi+R_2}}, \, \text{ and } \, CMY =\ln{\frac{\txtups+R_1}{\txtups+R_2}}.
\end{equation}
Moreover we can directly calculate 
$$CM5Z = \txtzet_1 \arctan{\frac{(\txtchi\txtups)_{ij}}{\txtzet_1(R_1)_{ij}}}, \, \text{ and }\,CM6Z=\txtzet_2\arctan{\frac{(\txtchi\txtups)_{ij}}{\txtzet_2(R_2)_{ij}}}.$$
Hence,  the summand of the triple sum is replaced by
\begin{align*}  \left(\left(  (CM5Z)_{ij} - (CM6Z)_{ij}\right)  -\txtchi_i (CMY)_{ij} -\txtups_j (CMX)_{ij}\right)
= CM_{ij}.
\end{align*}

Now
we see that $\txtchi_1$, $\txtchi_2$ are entries from $X$, $\txtups_1$, $\txtups_2$ are entries from $Y$. Thus, given the definition \eqref{PaddedDistances} and assuming that $X$ and $Y$ are stored in row vectors we can form matrices $XY=X(:).*Y$ and $R^2=X(:).^{\wedge }2+Y.^{\wedge }2$ which are of size $(\nbx+1) \times (\nby+1)$, and are independent of the $z$ coordinates. Thus, we save substantial computation by only calculating $XY$ and $R^2$ once for all slices, and for each slice we only calculate one row of the matrix.  Since these are based on matrices we can calculate the double sum for multiple coordinates by shifting each ${ij}$ matrix to the right in $i$ and right in $j$  with the appropriate sign, and obtain the entire sum in one line by correct indexing into the matrices. Suppose that $CM $ has size $(\nbx+1)\times (\nby+1)$ then we obtain a  matrix that can be reshaped to a row vector
\begin{align*}
g=CM(1:\nbx,1:\nby)-&CM(1:\nbx,2:\nby+1)-CM(2:\nbx+1,1:\nby)\\
+&CM(2:\nbx+1,2:\nbx+1)\\
\tilde{h}(a_1,b_1,0)=-\gamma &g(:).
\end{align*}
Thus, we calculate the first row of depth block $r$ by an evaluation of \eqref{htilde} for all coordinate contributions to station $1$ in one step. Note that the simplification \eqref{simplifiedCMX} is a further optimization of the calculation of entries for $G$ as compared to that given in \cite{ChenLiu:18}.
The details of the use and application of the of the gravity problem in the context of  forward algorithms for symBTTB matrices   are provided  in Algorithms~\ref{Alg:symBTTB}-\ref{Alg:symBTTBFFT} with the  $\gravity$ function in Algorithm~\ref{Alg:gravityresponse}. 
\subsection{Magnetic Kernel}\label{sec:magnetickernel}
We now extend the improvement of the calculation of the gravity kernel as discussed in  \cite{ChenLiu:18} to the calculation of the magnetic kernel, under the assumption that  there is no remanence magnetization or self-demagnetization, so that the magnetization vector is parallel to the  Earth's magnetic field\footnote{We assume that the total field is measured in nano Teslas;  introducing a scaling factor $10^9$ in the definitions.}.  Although the magnetic kernel is not symmetric, its discretization does lead to a BTTB matrix, and hence the discussion of Section~\ref{sec:BTTBstructure} is relevant. We apply the simplifications of  \cite{RaoBabu:91} for the evaluation of the magnetic kernel contribution. Using their notation \cite[eq. (3)]{RaoBabu:91}
\begin{equation}\label{magnetickernel}
\tilde{h}(a,b,0)_{pq}= \tilde{H}(G_1 \ln F_1 + G_2 \ln F_2 + G_3 \ln F_3 +G_4  F_4 +G_4 F_5).
\end{equation}
The constants $g_i=\tilde{H}G_i$ depend on the volume orientation and magnetic constants, \cite{RaoBabu:91}, and it is now assumed that $x$ points North and $y$ points East. Taking advantage of the notation in  \eqref{notationgravity}, the variables $F_i$ are given by 
\begin{align*}
F_1&=\frac{(\txtrho_{11}^2+\txtchi_1)(\txtrho_{21}^1+\txtchi_2)(\txtrho_{12}^1+\txtchi_1)(\txtrho_{22}^2+\txtchi_2)}
{(\txtrho_{11}^1+\txtchi_1)(\txtrho_{21}^2+\txtchi_2)(\txtrho_{12}^2+\txtchi_1)(\txtrho_{22}^1+\txtchi_2)}\\
F_2&=\frac{(\txtrho_{11}^2+\txtups_1)(\txtrho_{21}^1+\txtups_1)(\txtrho_{12}^1+\txtups_2)(\txtrho_{22}^2+\txtups_2)}
{(\txtrho_{11}^1+\txtups_1)(\txtrho_{21}^2+\txtups_1)(\txtrho_{12}^2+\txtups_2)(\txtrho_{22}^1+\txtups_2)}\\
F_3&=\frac{(\txtrho_{11}^2+\txtzet_2)(\txtrho_{21}^1+\txtzet_1)(\txtrho_{12}^1+\txtzet_1)(\txtrho_{22}^2+\txtzet_2)}
{(\txtrho_{11}^1+\txtzet_2)(\txtrho_{21}^2+\txtzet_1)(\txtrho_{12}^2+\txtzet_1)(\txtrho_{22}^1+\txtzet_2)}\\
F_4=\arctan\frac{\txtchi_2\txtzet_2}{\txtrho_{22}^2\txtups_2}&-\arctan\frac{\txtchi_1\txtzet_2}{\txtrho_{12}^2\txtups_2}-\arctan\frac{\txtchi_2\txtzet_2}{\txtrho_{21}^2\txtups_1}+
\arctan\frac{\txtchi_1\txtzet_2}{\txtrho_{11}^2\txtups_1}-\\
&
\arctan\frac{\txtchi_2\txtzet_1}{\txtrho_{22}^1\txtups_2}+
\arctan\frac{\txtchi_1\txtzet_1}{\txtrho_{12}^1\txtups_2}+
\arctan\frac{\txtchi_2\txtzet_1}{\txtrho_{21}^1\txtups_1}-
\arctan\frac{\txtchi_1\txtzet_1}{\txtrho_{11}^1\txtups_1}\\
F_5=\arctan\frac{\txtups_2\txtzet_2}{\txtrho_{22}^2\txtchi_2}&-\arctan\frac{\txtups_2\txtzet_2}{\txtrho_{12}^2\txtchi_1}-\arctan\frac{\txtups_1\txtzet_2}{\txtrho_{21}^2\txtchi_2}+
\arctan\frac{\txtups_1\txtzet_2}{\txtrho_{11}^2\txtchi_1}-\\
&
\arctan\frac{\txtups_2\txtzet_1}{\txtrho_{22}^1\txtchi_2}+
\arctan\frac{\txtups_2\txtzet_1}{\txtrho_{12}^1\txtchi_1}+
\arctan\frac{\txtups_1\txtzet_1}{\txtrho_{21}^1\txtchi_2}-
\arctan\frac{\txtups_1\txtzet_1}{\txtrho_{11}^1\txtchi_1}.
\end{align*}
As for the gravity kernel, the calculation of \eqref{magnetickernel} can, therefore, use \eqref{notationgravity} to calculate $X$, $Y$ and $R^2$ once for all slices. We note that minor computational savings may be made by calculating for example $R_1+\txtchi_1$ within the calculations for $\tilde{h}$ but these are not calculations that can be made independent of the given slice. It may appear also that one could calculate ratios $\txtchi/\txtups$, with modification of the calculations for $F_4$ and $F_5$, but the stable calculation of the $\arctan$ requires the ratios as given. Otherwise sign changes in the numerator or denominator passed to $\arctan$ can lead to changes in the obtained angle. In MATLAB we use \texttt{atan2} rather than \texttt{atan} for improved stability in the calculation of the angle.  The details of the use and application of the of the magnetic problem in the context of the forward algorithms for  BTTB matrices   are provided in Algorithms~\ref{Alg:BTTB}-\ref{Alg:BTTBFFT}  with the  $\magnetic$  function in Algorithm~\ref{Alg:magneticresponse}.

\section{Numerical Validation}\label{sec:numerics}

We now validate the fast and efficient methods for generating both the symmetric and non symmetric kernels relating to gravity and magnetic problems.  Specifically, we compare the computational cost of direct calculation of the entries of the matrix $G$ that are required for matrix multiplications, with the entries that are required for the transform implementation of the multiplications. Thus, we compare Algorithms~\ref{Alg:symBTTB} and \ref{Alg:symBTTBFFT} with all entries calculated using Algorithm~\ref{Alg:gravityresponse}, and Algorithms~\ref{Alg:BTTB} and \ref{Alg:BTTBFFT} with all entries calculated using Algorithm~\ref{Alg:magneticresponse},  for the symmetric gravity, and non symmetric magnetic kernels, respectively. The $12$ problem sizes considered are detailed in Table~\ref{Table1}. They are generated by taking the smallest size with $(\nsx,\nsy,,\nbz) = (25,15,2)$, and then scaling each dimension by $1$ to $12$ for the test cases. For padding we compare the case with $p=0\%$ and $p=5\%$ padding across $x$ and $y$ dimensions, and rounded to the nearest integer. Thus, $m = \nsx\nsy$, and $n = \lfloor(1+p)\nsx\rfloor\lfloor(1+p)\nsy\rfloor\nbz$. All computations use MATLAB release  2019b implemented on a desktop computer with  an Intel(R) Xeon (R) Gold $6138 $ processor  ($2.00$GHz) and  $256$GB RAM.
\begin{table}[htb!]\begin{center}
\begin{tabular}{|c|c|c|c|c|}
\hline
Problem& $(\nsx,\nsy,,\nbz)$ & $m$ & $n$ ($p = 0\%$) & $n$ ($p = 0.05\%$)\\
\hline
$1$& $(25,   15,    2)  $&$ 375     $&$ 750         $&$ 918        $\\ \hline
$2$& $(50,   30,    4)  $&$ 1500   $&$ 6000       $&$ 7616      $\\ \hline
$3$& $(75,   45,    6)  $&$ 3375   $&$ 20250     $&$ 24402    $\\ \hline
$4$& $(100, 60,    8)  $&$ 6000   $&$ 48000     $&$ 58080    $\\ \hline
$5$& $(125, 75,   10) $&$ 9375   $&$ 93750     $&$ 113710  $\\ \hline
$6$& $(150, 90,   12) $&$ 13500 $&$ 162000   $&$ 199200  $\\ \hline
$7$& $(175, 105, 14) $&$ 18375 $&$ 257250   $&$ 310730  $\\ \hline
$8$& $(200, 120, 16) $&$ 24000 $&$ 384000   $&$ 464640  $\\ \hline
$9$& $(225, 135, 18) $&$ 30375 $&$ 546750   $&$ 662450  $\\ \hline
$10$& $(250, 150, 20) $&$ 37500 $&$ 750000   $&$ 916320  $\\ \hline
$11$& $(275, 165, 22) $&$ 45375 $&$ 998250   $&$ 1206500$\\ \hline
$12$& $(300, 180, 24) $&$ 54000 $&$ 1296000 $&$ 1568200$\\ \hline
\end{tabular}
\caption{Dimensions of the volume used in the experiments labeled as problems $1$ to $12$ corresponding to scaling each dimension in $(25,   15,    2)  $ by the problem number and increasing $m$ by a factor $8$ for each row. \label{Table1}}
\end{center}\end{table}

In the results, we reference the kernels generated by Algorithms~\ref{Alg:symBTTB}, \ref{Alg:BTTB}, \ref{Alg:symBTTBFFT}, and \ref{Alg:BTTBFFT} as $\Ggrav$, $\Gmag$, $\Tgrav$, and $\Tmag$ respectively. These values are plotted on a ``log-log" scale  in Figure~\ref{figure5}, without padding in Figure~\ref{figure5a} and with padding in Figure~\ref{figure5b},  in which the problem sizes are given as relevant triples on the $x-$axis. The   problem cases  from $8$ to $12$ for the direct calculation of $G$ are too large to fit in memory on the given computer.  It can be seen that the generation of $G$ is effectively independent of the \gravity ~or \magnetic~kernels; $\Ggrav$, $\Gmag$ are comparable. But the requirement to calculate extra entries for the non-symmetric \magnetic~kernel is also seen; $\Tgrav<\Tmag$. On the other hand, the significant savings in generating just the transform matrices, as indicated by timings $\Tgrav$, and $\Tmag$, as compared to $\Ggrav$, and $\Gmag$ is evident. There is a considerable computational advantage to the use of the transform for calculating the required components that are needed for evaluating matrix-vector products for these structured kernel matrices. 

\begin{figure}[ht!]\begin{center}
\subfigure[$p=0$ \label{figure5a}]{\includegraphics[width=0.49\textwidth]{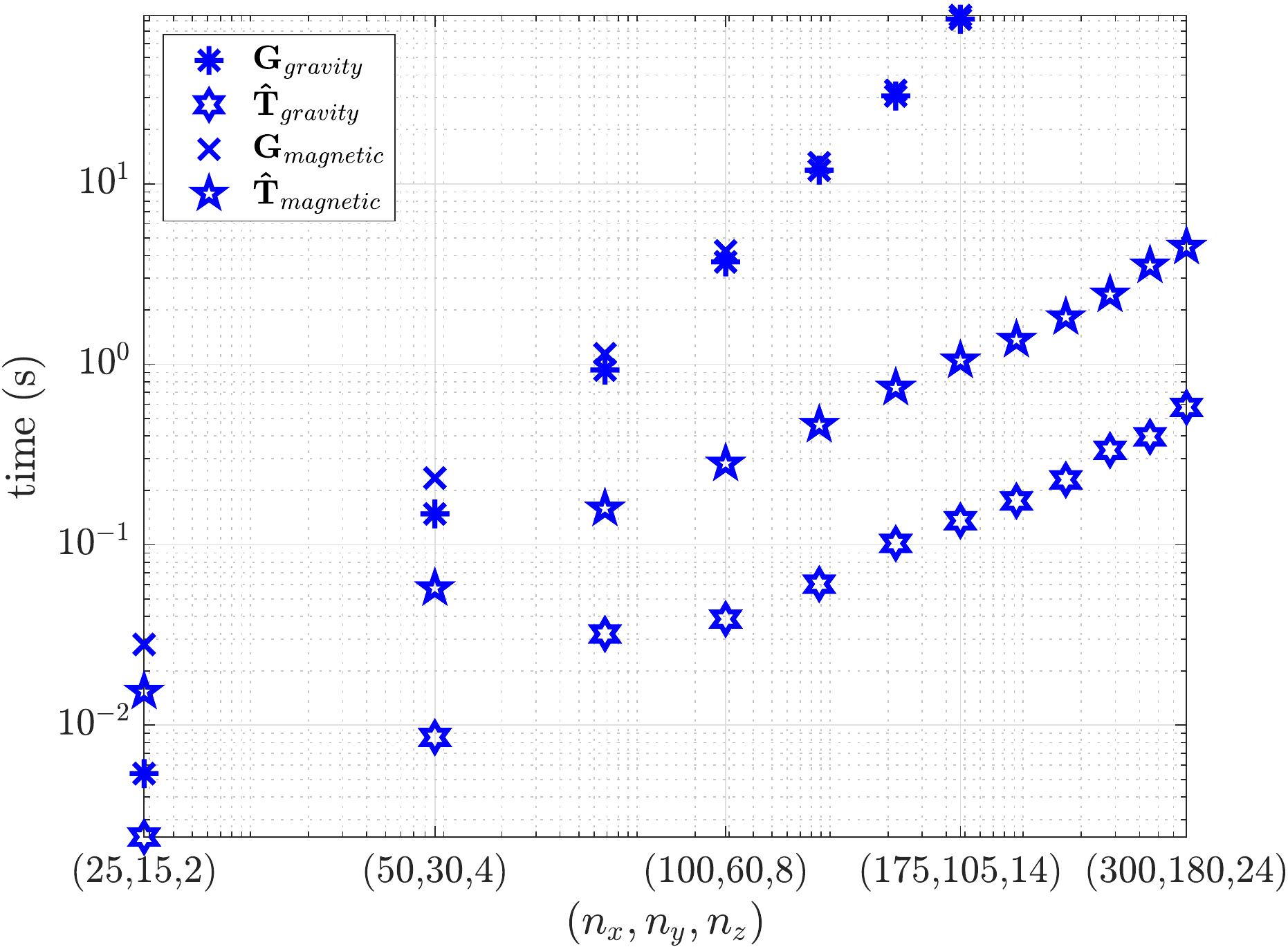}}
\subfigure[$p=0.05$ \label{figure5b}]{\includegraphics[width=0.49\textwidth]{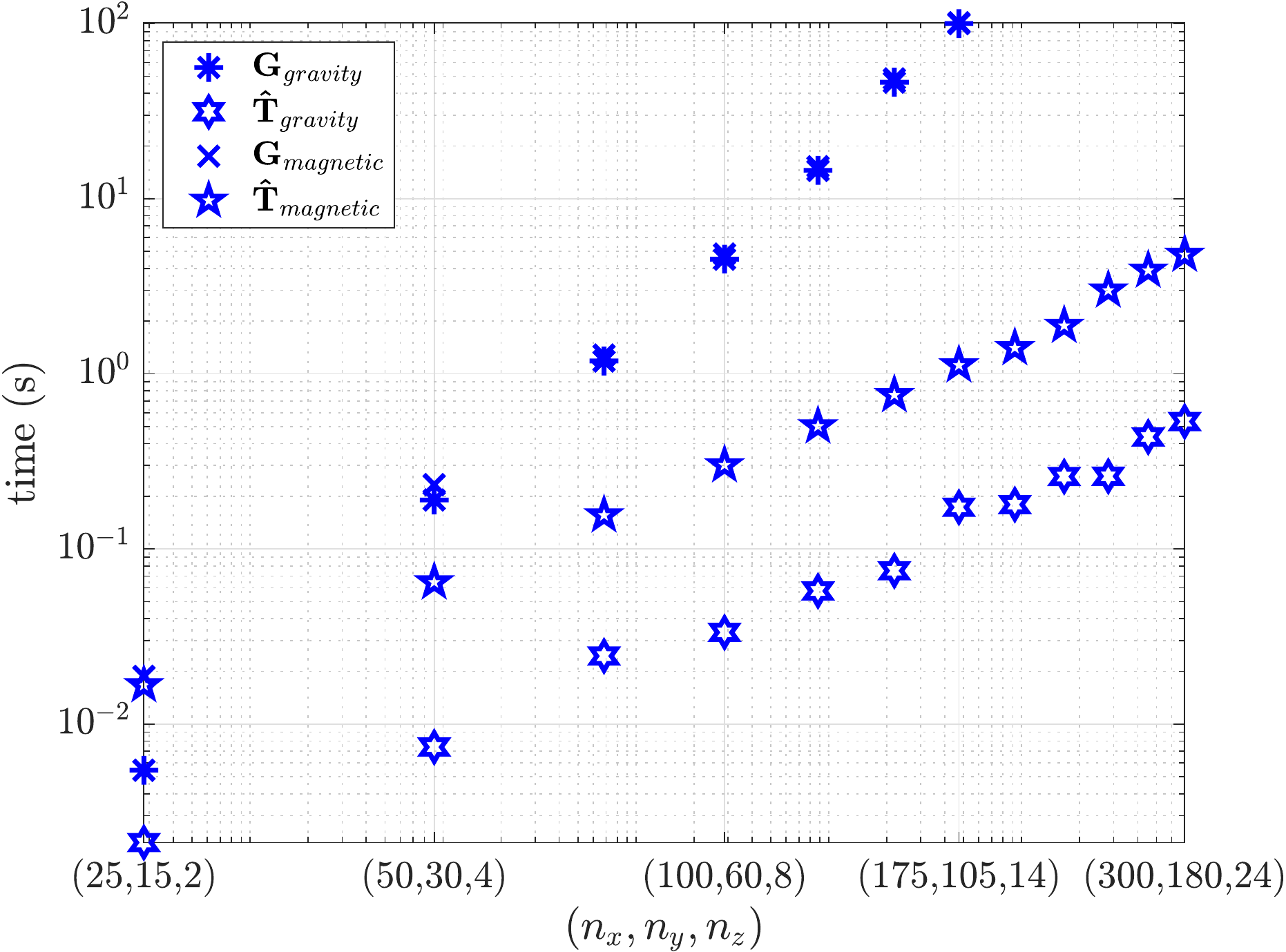}}
\caption{Running times   for generating $G$ using Algorithms~\ref{Alg:symBTTB} and \ref{Alg:BTTB} and $\hat{T}$ using Algorithms~\ref{Alg:symBTTBFFT} and \ref{Alg:BTTBFFT} with $0\%$ padding in Figure~\ref{figure5a}, and with $5\%$ padding in Figure~\ref{figure5b}.\label{figure5}}
\end{center} \end{figure}

Of greater significance  is the comparison of the computational cost of  direct matrix multiplications, $\bfb=G \bfu$ and $\bfd=G^T \bfv$,  as compared with the transform implementations for these products. 
As noted already in \eqref{Gslice}, we assume that the matrices decompose by depth  so that forward and transpose operations are implemented by depth level, as indicated in Algorithm~\ref{Alg:multbttb}. Thus, we  consistently partition   $\bfu \in \Rm{\nbx\nby\nbz}$ into $\nbz$ blocks,    $\bfu_r\in\Rm{\nbx\nby}$, $1\le r\le \nbz$. Then, 
\begin{equation*}
G\bfu = \sum_{r = 1}^{\nbz}\Gr\bfu_r, \text{ and }\\
G^T\bfv = \left[\begin{array}{c} \left(G^{(1)}\right)^T\bfv\\
\left(G^{(2)}\right)^T\bfv\\
\vdots\\
\left(G^{(\nbz)}\right)^T\bfv\end{array}\right],
\end{equation*}
where $\bfv\in\Rm{\nsx\nsy}$. 
$100$ copies of vectors $\bfu\in\Rm{n}$ and $\bfv\in\Rm{m}$ are randomly generated and the mean times for calculating the products over all    $100$ trials, for each problem size, are recorded.
We also record the differences over all trials in the generation of  $\bfb$ and $\bfd$ obtained directly for $\Ggrav$ and $\Gmag$ and by Algorithm~\ref{Alg:multbttb} for $\Tgrav$ and $\Tmag$.  Then,  $E_{\gravity}$ and $E_{\magnetic}$ are the mean values of  the relative $2$-norm of the difference between the results produced by $\Ggrav$ versus $\Tgrav$, and for $\Gmag$ versus $\Tmag$, respectively, for both forward and transpose operations. The results are illustrated in Figures~\ref{figure6} and  \ref{figure7} for the generation of $G\bfu$ and $G^T\bfv$, respectively. In each case the timing is reported on the left $y-$axis and the error on the right $y$-axis. Again all plots are on the ``log-log" scale, and Figures~\ref{figure6a} and  \ref{figure7a} are without padding, but Figures~\ref{figure6b} and  \ref{figure7b} are with padding. Figures~\ref{figure6} and Figure~\ref{figure7} show significant reductions in mean running time when implemented without the direct calculation of the matrices. Moreover, the results are comparable,    $E_{\gravity} \lesssim 10 \epsilon$ for both forward and transpose operations, and $E_{\magnetic} \lesssim 10^2 \epsilon$, where $\epsilon$ is the machine accuracy.      Thus, in all cases, $\Tgrav$ and $\Tmag$ show a significant reduction in mean running time for large problems, and allow much larger systems to be represented. Indeed, the largest test case for $\Tgrav$ and $\Tmag$ is by no means a limiting factor, and it is possible to represent much larger kernels.

\begin{remark}[MATLAB $\fftt$]
We should note that the MATLAB $\fftt$ function determines an optimal algorithm for a given problem size. On the first call for a given problem size, $\fftt$ uses the function $\fftw$ to determine optimal parameters for the Fourier transform.  Thus, the first time $\fftt$ is called generally takes longer than subsequent instances.  We mitigate this effect by first removing the variable $\dwisdom$ within $\fftw$, and then setting the planner within $\texttt{fftw}$ to exhaustive.  After running a single call to $\fftt$, the resulting $\dwisdom$ is then saved.  This process is repeated for generating $\hat{T}$, forward multiplication, and transpose multiplication.  Then for each trial, the appropriate stored values for $\dwisdom$ are loaded before each use of $\fftt$. Hence the results are not contaminated by artificially high costs of the first run of $\fftw$ for each problem case. 
\end{remark}

\begin{figure}[ht!]\begin{center}
\subfigure[$p=0$ \label{figure6a}]{\includegraphics[width=0.36\textwidth,angle=90]{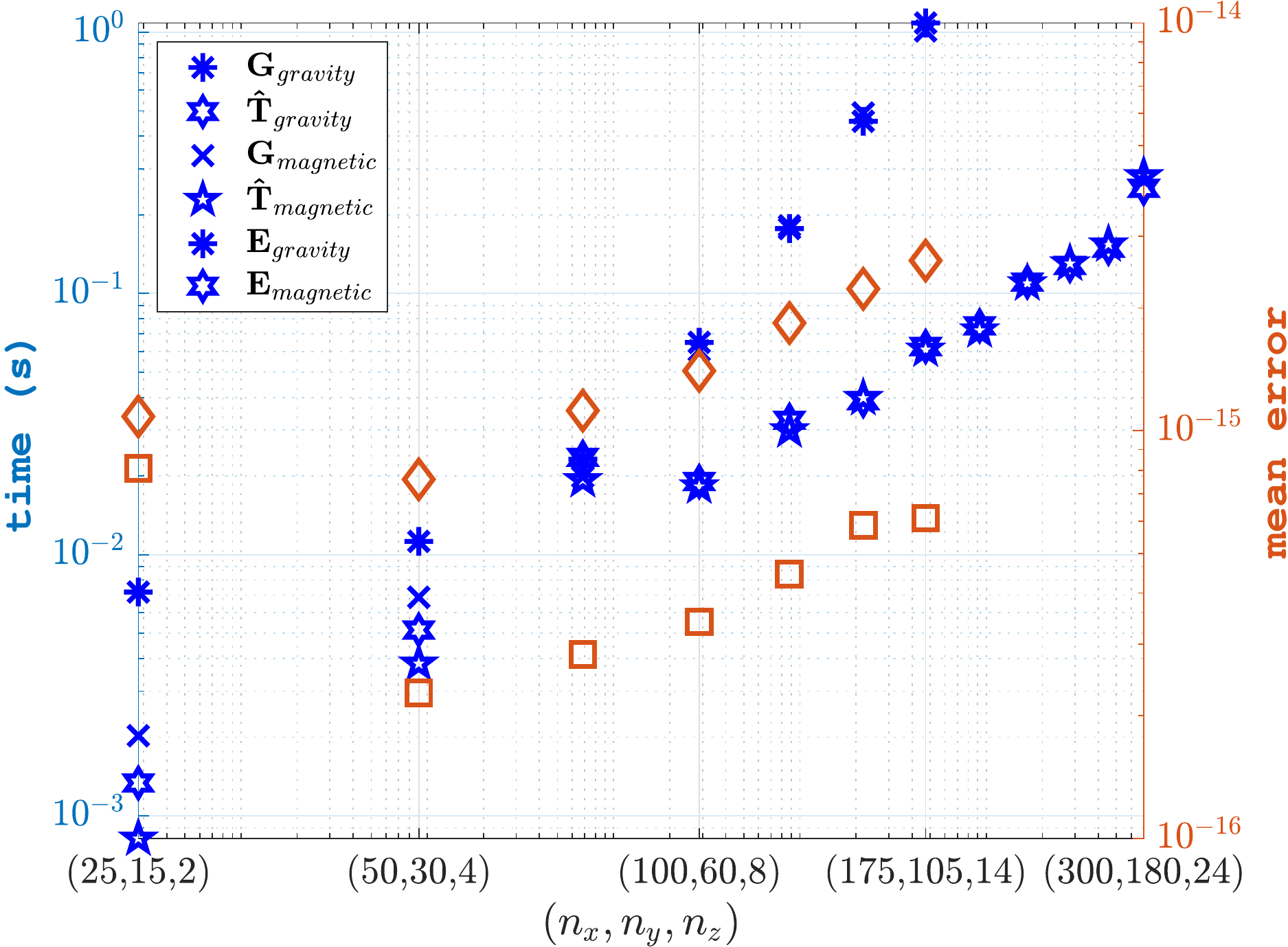}}
\subfigure[$p=0.05$ \label{figure6b}]{\includegraphics[width=0.36\textwidth,angle=90]{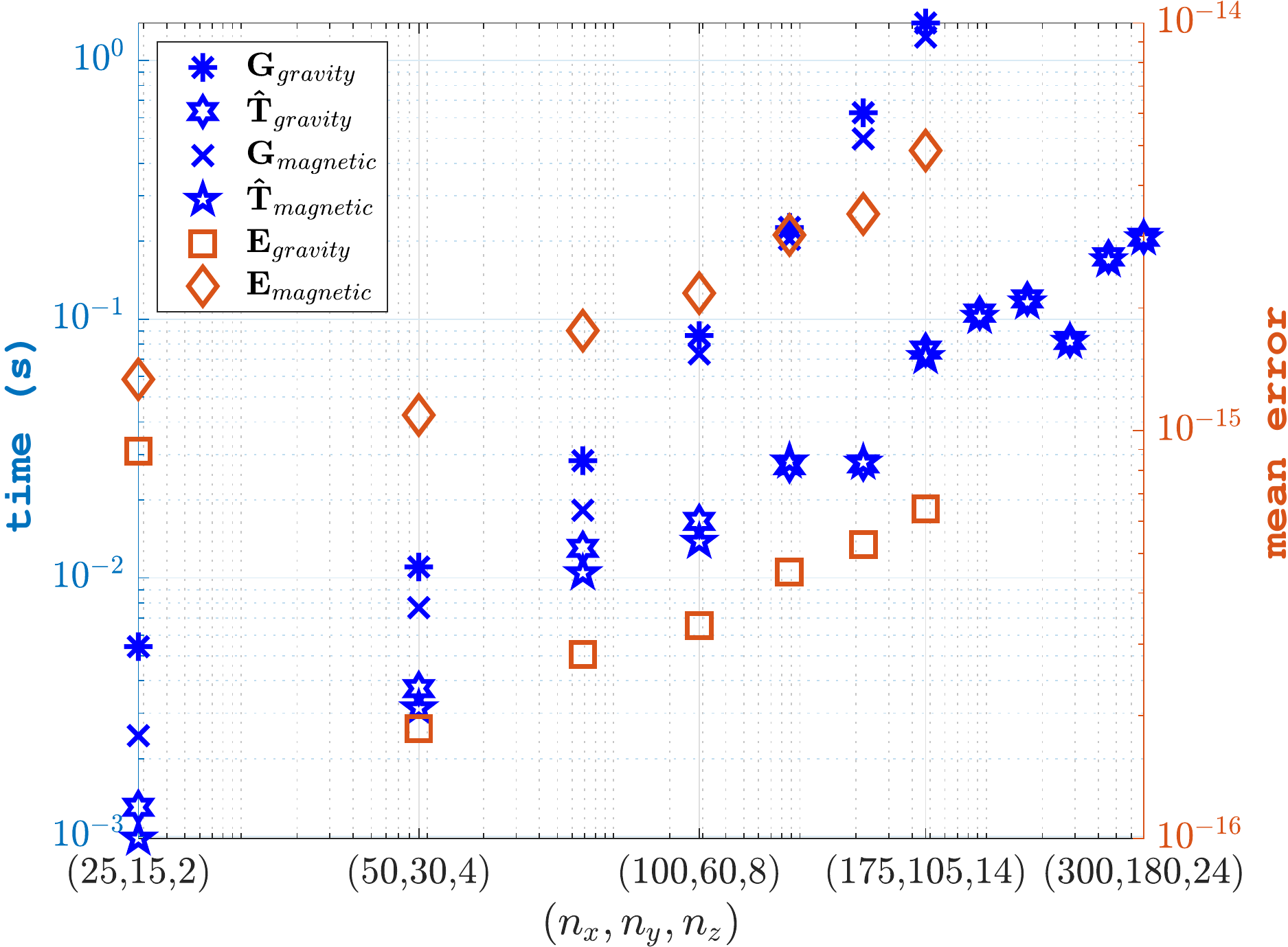}}
\caption{Mean running times over $100$ trials for forward m ultiplication using $\Ggrav$, $\Gmag$, $\Tgrav$, and $\Tmag$ (left $y$-axis) and mean errors over $100$ trials $E_{\gravity}$ and $E_{\magnetic}$ for forward multiplication using $\Ggrav$ versus $\Tgrav$, and $\Gmag$ versus $\Tmag$ respectively (right $y$-axis) are shown for $0\%$ padding in Figure~\ref{figure6a}, and $5\%$ padding in Figure~\ref{figure6b}.\label{figure6}}
\end{center}\end{figure}

\begin{figure}[ht!]\begin{center}
\subfigure[$p=0$ \label{figure7a}]{\includegraphics[width=0.36\textwidth,angle=90]{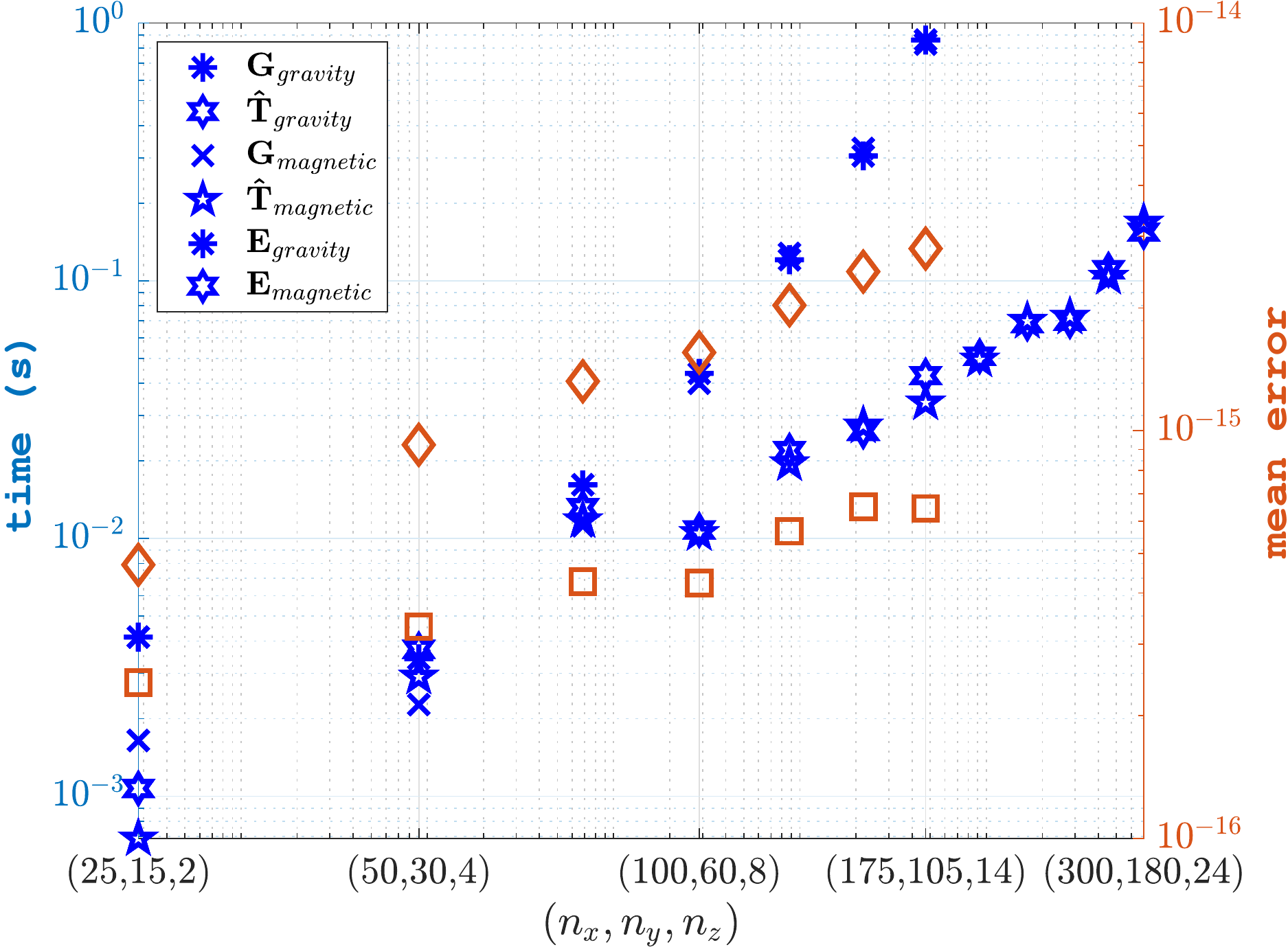}}
\subfigure[$p=0.05$ \label{figure7b}]{\includegraphics[width=0.36\textwidth,angle=90]{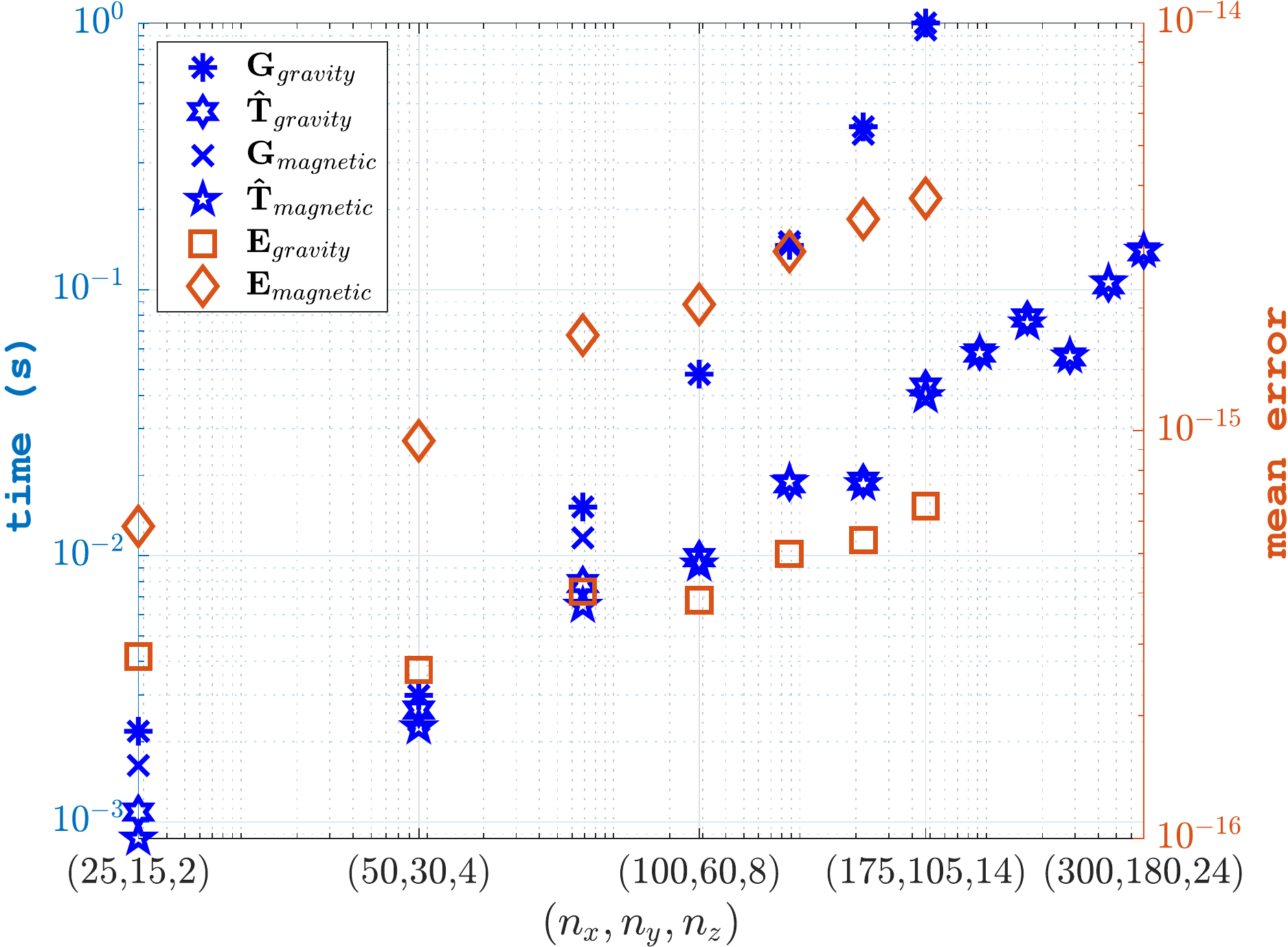}}
\caption{Mean running times over $100$ trials for transpose multiplication using $\Ggrav$, $\Gmag$, $\Tgrav$, and $\Tmag$ (left $y$-axis) and mean errors over $100$ trials $E_{\gravity}$ and $E_{\magnetic}$ for transpose multiplication using $\Ggrav$ versus $\Tgrav$, and $\Gmag$ versus $\Tmag$ respectively (right $y$-axis) are shown for $0\%$ padding in Figure~\ref{figure7a}, and $5\%$ padding in Figure~\ref{figure7b}.\label{figure7}}
\end{center}\end{figure}

\section{Data Availability and Software Package}\label{sec:code}
The software consists of the main functions to calculate the BTTB and symBTTB matrices, with padding, and the circulant matrices $\Tcirc$ that are needed for the 2DFFT.  Also provided is a simple script to test the algorithms using the gravity and magnetic kernels.  All the algorithms are described in Appendix~\ref{App:B} and the software is  is open source and available at  \url{https://github.com/renautra/FastBTTB}, and described at  \url{https://math.la.asu.edu/~rosie/research/bttb.html}. Provided are the scripts that are used to generate the results presented in the paper. The variables used in the codes are described in Table~\ref{tab:notation} and Table~\ref{tab:parameters}. The \texttt{TestingScript}.m is easily modified to generate new examples and can be tested within different hardware configurations and versions of Matlab. A safety test for memory usage in generating large scale examples is provided at the initialization of each problem size, so that problems too large to fit in memory will not be used in generating the matrix $G$ directly.  The presented implementation assumes uniform grid sizes in all dimensions, but using depth layers of different heights is an easy modification, through the change in the input coordinate vector $z\_blocks$.

\section{Conclusions and Future Work}\label{sec:conclusions}
We have provided  a description of the generation of efficient codes for implementing forward and transpose operations with BTTB matrices. These are used in geophysical forward modeling for cases in which the kernels are of convolution type and generate matrices with the required structure. Then, efficient generation of matrix operations with minimal storage makes it feasible to perform large three-dimensional modeling with these kernels. A novelty of this work, beyond existing descriptions in the literature, is the generation of operators that allow for padding in the coordinate volume. Also, in this work, the approach for finding operations $G^T\bfm$ explicitly given knowledge of the BTTB structure of $G$ and its BCCB embedding is provided. Thus, the developed software can be integrated into an inverse modeling problem, in which given data $\bfd$,  model parameters $\bfm$ are desired. This is planned for future work.  

\section*{Acknowledgments}
Rosemary Renaut acknowledges the support of NSF grant  DMS 1913136:   ``Approximate Singular Value Expansions and Solutions of
Ill-Posed Problems".

\appendix
\section{Notation and Parameter Definitions}\label{App:A}
\begin{table}[htb!]
\begin{center}
\begin{tabular}{|c|c |}\hline
$\param$& $\nsx,\nsy,\nbz,\padxl,\padxr,\padyl,\padyr,\nbx,\nby,m,n,\nbr,\padx,\pady$ \\ \hline
$gsx$, $gsy$, $gsz$ & Grid sizes $\dx$, $\dy$ and $\dz$   \\ \hline
$\That$ & $\Thatforward=\hat{T}^{\mathrm{circ}}$, $\Thattranspose=\hat{\tilde{T}}^{\mathrm{circ}}$ \\ \hline
z$\_$blocks & Depth coordinates, increasing, $z_r$ \\ \hline
$D$ & Declination of geomagnetic field and  magnetization vector  \\ \hline
  $I$ & Inclination of  geomagnetic field and  magnetization vector  \\ \hline
$F$ &  Intensity of the geomagnetic field  in $\mathrm{nT}$ ($10^{-9}F$ in $\mathrm{T}$)\\ \hline
$H=\frac{10^{-9}F}{4\pi}$ & Magnetic field intensity ($\mathrm{A}/\mathrm{m}$) in SI units \\ \hline
$\tilde{H}=10^9 H =\frac{F}{4\pi}$ & Assumes the field is measured in $\mathrm{nT}$  \\ \hline
\end{tabular}
\caption{Parameters and variables  in the codes. The parameters are defined in Table~\ref{tab:notation}\label{tab:parameters}}
\end{center}
\end{table}

\begin{table}[htb!]
\begin{center}
\begin{tabular}{|c|c|c|c|}
\hline
$\nsx$& $\#$ true stations in $x$&$s_{ij}=(a_{ij},b_{ij})$& Station location\\ \hline 
$\nsy$& $\#$ true stations in $y$ &$m=\nsx\nsy$& $\#$ measurements\\\hline
$\nbx,\nby,\nbz$&$\#$ coordinate blocks in $x$,$y$,$z$&
$\dx,\dy,\dz$& Grid sizes in $x$,$y$,$z$ \\ \hline
$\padxl,\padxr,\padx$& Left, right, total padding: $x$&$\nbx$& $\nbx=\nsx+\padx$ \\ \hline 
$\padyl,\padyr,\pady$& Left, right,total, padding: $y$&$\nby$& $\nby=\nsy+\pady$\\ \hline
$x_p$& $x_p=(p-1-\padxl)\dx$, $1\le p \le \nbx+1$  &$n=\nbx\nby\nbz$& Volume Dimension \\ \hline
$y_q$& $y_q=(q-1-\padyl)\dy$,  $1\le q \le \nby+1$& $\nbr=\nbx\nby$& Layer Dimension \\ \hline 
$z_r$&   $z_r=(r-1)\dz$, $ 1\le r \le \nbz+1$&$c_{pqr}$&Prism  $pqr$ in $xyz$ \\ \hline
 $d$, $h$, $\zeta$&Forward Model  see \eqref{continuousforwardmodel}&
$\tilde{h}(s_{ij})_{pqr}$& Projection  $c_{pqr}$ to   $s_{ij}$\\ \hline
$G\in\Rmn{m}{n} $& $(\Gr)_{k\ell}=\tilde{h}(s_{ij})_{pqr}$ See \eqref{htilde}&   $\Gr\in\Rmn{m}{\nbr} $& 
Depth $r$  Contribution\\ \hline 
  $\Gr_q\in\Rmn{\nsx}{\nbx}$&  $\Gr_q=\Gr_{1q}$, $1\le q \le \nby$ &$\barGr_j\in\Rmn{\nsx}{\nbx}$& $\barGr_j=\barGr_{j1}$, $1\le j \le \nsy$ \\ \hline
 $\bfc_q$, $\bfr_q$&$\Gr_q=\mathrm{toeplitz}(\bfc_q,\bfr_q)$& $\barc_j$, $\barr_j$&  $\barGr_j=\mathrm{toeplitz}(\barc_j,\barr_j)$ \\ \hline
  BTTB &Block Toeplitz Toeplitz blocks &symBTTB & Symmetric BTTB\\ \hline 
    BCCB &Block Circulant Circulant blocks &   $J_m$ Defn.~\ref{defn:exchange}& Exchange matrix\\ \hline
$\bcext_q$, $\brext_q$& Defining $(\Gr)^{\mathrm{circ}}$& $\barcext_j$, $\barrext_j$&  Defining $(\barGr)^{\mathrm{circ}}$\\ \hline
    $\Tcirc$& Components of BCCB &  $\tilde{T}^{\mathrm{circ}}$& Transpose Components   \\ \hline 
    $\hat{T}^{\mathrm{circ}}$ & $\fftt(\Tcirc)$ : 2DFFT&  $\hat{\tilde{T}}^{\mathrm{circ}}$ & $\fftt(\tilde{T}^{\mathrm{circ}})$ : 2DFFT\\ \hline

\end{tabular}
\end{center}
\caption{Notation Adopted in the Discussion \label{tab:notation}}
\end{table}
\section{Algorithms}\label{App:B}
An overview of the  required algorithms as described in Sections~\ref{sec:forwardmodel}-\ref{sec:circulant} are provided in Algorithms~\ref{Alg:symBTTB}-\ref{Alg:symBTTBFFT} using the  $\gravity$  function in Algorithm~\ref{Alg:gravityresponse} and in Algorithms~\ref{Alg:BTTB}-\ref{Alg:BTTBFFT}  with the  $\magnetic$  function in Algorithm~\ref{Alg:magneticresponse}. The convolution multiplication is provided in Algorithm~\ref{Alg:multbttb}. 

\begin{center}
\begin{algorithm}[H] 
\caption{$G=\symbttb(gsx,gsy,z\_blocks,\param)$\newline\label{Alg:symBTTB}
Entries of padded symBTTB matrix. Function \gravity.}
\SetAlgoLined \LinesNumbered
  \KwIn{See Table~\ref{tab:parameters} for details\;
  $gsx$, $gsy$, $z\_blocks$ : grid spacing in $x$ and $y$ and $z$ coordinates\;
  $\param$ required parameters } 
\KwOut {\
symBTTB real matrix $G$ of size $m \times n$.}
Extract parameters from $\param$ \;
Initialize zero arrays: $G$, $Gr$, $Grq$\;
Sizes: $nX=\nsx+\max(\padxl,\padxr)$, $nY=\nsy+\max(\padyl,\padyr)$\;
Form distance arrays $X$ and $Y$ according to \eqref{shortsymXY}\;
Form  $X2=X^2$, $Y2=Y^2$, $XY=X(:).*Y$ and $R=X2(:)+Y2$\;
\For{$r=1:\nbz$ }{
Set $z_1=z\_blocks(r)$, $z_2=z\_blocks(r+1)$\;
Calculate slice response at first station: $g=\gravity(z_1,z_2, X,Y,XY,R)$\;
\For{$q=1:nY$}{
            Extract $\bfc_q$, $\bfr_q$ from $g$ :  use \eqref{cqpaddedsymtoeplitz}, \eqref{rqpaddedsymtoeplitz}\;
            Generate: $\Gr_q=\mathrm{toeplitz}(\bfc_q,\bfr_q)$: use \eqref{paddedsymtoeplitz} \;}
\For{$j=[\padyl+1:-1:2,1:\nsy+\padyr]$}{
	    Build row first row of $Gr$ using \eqref{Rpaddedsymtoeplitz}\;}
\For{$j=2:\nsy$}{ 
            Build $j^{\mathrm{th}}$ row of $\Gr$ using \eqref{paddedsymmetrictoeplitzblockmatrix} and \eqref{Cpaddedsymtoeplitz}\;}
Assign: $Gr$ to $r^{\mathrm{th}}$ block of $G$\;}
\end{algorithm}
\end{center}

\begin{center}
\begin{algorithm}[htb!]
\caption{$\That=\symbttbfft(gsx,gsy,z\_blocks,\param)$\newline\label{Alg:symBTTBFFT}
Transforms of padded symBTTB matrix. Function \gravity.}
\SetAlgoLined \LinesNumbered
  \KwIn{See Table~\ref{tab:parameters} for details\;
  $gsx$, $gsy$, $z\_blocks$ : grid spacing in $x$ and $y$ and $z$ coordinates\;
  $\param$ required parameters } 
\KwOut {\
Structure $\That$ with $\Thatforward$ and $\Thattranspose$}
Extract parameters from $\param$ \;
Initialize zero arrays: $\Tcirc$, $\hat{T}^{\mathrm{circ}}$ and $\hat{\tilde{T}}^{\mathrm{circ}}$\;
Sizes: $nX=\nsx+\max(\padxl,\padxr)$, $nY=\nsy+\max(\padyl,\padyr)$\;
Form distance arrays $X$ and $Y$ according to \eqref{shortsymXY}\;
Form  $X2=X^2$, $Y2=Y^2$, $XY=X(:).*Y$ and $R=X2(:)+Y2$\;
\For{$r=1:\nbz$ }{
Set $z_1=z\_blocks(r)$, $z_2=z\_blocks(r+1)$\;
Calculate slice response at first station: $g=\gravity(z_1,z_2, X,Y,XY,R)$\;
\For{$j=[1+\padyl:\nsy+\padyl. \nsy+\padyr:-1:2 1:\padyl]$}{
            Extract $\bfr_j$ from $g$ :  use  \eqref{rqpaddedsymtoeplitz}\;
            Augment column of  $\Tcirc$, use \eqref{Eqn:Tcircpadded} \;}
Take FFT of  $\Tcirc$: $\hat{T}^{\mathrm{circ}}(:,:,r)=\fftt(\Tcirc)$ \;
Generate $T$ from $\Tcirc$ using \eqref{Eqn:Ttran0} and \eqref{Eqn:Ttran}\;
Take FFT of  $T$: $\hat{\tilde{T}}^{\mathrm{circ}}(:,:,r)=\fftt(T)$ \;}
$\Thatforward=\hat{T}^{\mathrm{circ}}$, $\Thattranspose=\hat{\tilde{T}}^{\mathrm{circ}}$.
\end{algorithm}
\end{center}

\begin{center}
\begin{algorithm}[htb!] 
\caption{$G=\bttb(gsx,gsy,z\_blocks,\param, D, I, H)$\newline\label{Alg:BTTB}
Entries of padded BTTB matrix,  Figure~\ref{figure3}. Function \magnetic.}
\SetAlgoLined \LinesNumbered
  \KwIn{See Table~\ref{tab:parameters} for details\;
  $gsx$, $gsy$, $z\_blocks$ : grid spacing in $x$ and $y$ and $z$ coordinates\;
  $\param$ required parameters } 
\KwOut {\
BTTB real matrix $G$ of size $m \times n$. }
Extract parameters from $\param$ \;
Calculate constants $g_i=\tilde{H}G_i$ for \eqref{magnetickernel}, \cite[(3)]{RaoBabu:91}\;
Initialize zero arrays: $G$, $Gr$ and row and column cell arrays,   Figure~\ref{figure3}\;
Sizes: $nX=\nsx+\max(\padxl,\padxr)$, $nY=\nsy+\max(\padyl,\padyr)$\;
Form distance arrays $X$ and $Y$ according to \eqref{shortunsymXY}\;
Form  $X2=X^2$, $Y2=Y^2$, and $R=X2(:)+Y2$\;
\For{$r=1:\nbz$ }{
Set $z_1=z\_blocks(r)$, $z_2=z\_blocks(r+1)$\;
Calculate $grow\{1\} = \tilde{h}({11})_{pq}$, $1\le p \le nX$, $1\le q\le nY$\;
\For{$j=2:\nsy+\padyl$}{
Calculate $grow\{j\} = \tilde{h}({1j})_{p1}$, $1\le p \le nX$\;
}
Calculate : $gcol\{1\}= \tilde{h}({ij})_{11}$, $1\le i \le nX$, $1\le j\le nY$\;
\For{$q=2:\nsy+\padyr$}{
Calculate : $gcol\{q\}= \tilde{h}({i1})_{1q}$, $1\le i \le nX$\;}
\For{$j=\padyl+1:-1:2$}{
	    Generate  $\barGr_j$: using $gcol\{1\}$ and $grow\{j\}$, \eqref{barpaddedtoeplitz} for $R$ in \eqref{nonsymBTTBtoeplitzpadded}\;}
\For{$q=1:\nsy+\padyr$}{
	    Generate $\Gr_q$: using $gcol\{q\}$ and $grow\{1\}$, \eqref{paddedtoeplitz} for  $R$ in \eqref{nonsymBTTBtoeplitzpadded}\;}
\For{$j=\padyl+2:-1:2,1:\nsy+\padyl$}{
	    Generate  $\barGr_j$: using $gcol\{1\}$ and $grow\{j\}$, \eqref{barpaddedtoeplitz} for $C$ in \eqref{nonsymBTTBtoeplitzpadded}\;}
Build $\Gr$  in \eqref{nonsymBTTBtoeplitzpadded} using $C$ and $R$\;
Assign: $Gr$ to $r^{\mathrm{th}}$ block of $G$\;}
\end{algorithm}
\end{center}

\begin{center}
\begin{algorithm}[htb!] 
\caption{$\That=\bttbfft(gsx,gsy,z\_blocks,\param, D, I, H)$\newline
Transforms of padded BTTB matrix,  Figure~\ref{figure3}. Function \magnetic.\label{Alg:BTTBFFT}}
\SetAlgoLined \LinesNumbered
  \KwIn{See Table~\ref{tab:parameters} for details\;
  $gsx$, $gsy$, $z\_blocks$ : grid spacing in $x$ and $y$ and $z$ coordinates\;
  $\param$ required parameters\; $D$, $I$, $H$ declination, inclination and intensity of magnetization} 
\KwOut {\
Structure $\That$ with $\Thatforward$ and $\Thattranspose$}
Extract parameters from $\param$ \;
Calculate constants $g_i=\tilde{H}G_i$ for \eqref{magnetickernel}, \cite[(3)]{RaoBabu:91}\;
Initialize zero arrays: $\Tcirc$, $\hat{T}^{\mathrm{circ}}$ and $\hat{\tilde{T}}^{\mathrm{circ}}$\;
Initialize zero arrays for and row and column cell arrays, see Figure~\ref{figure3}\;
Sizes: $nX=\nsx+\max(\padxl,\padxr)$, $nY=\nsy+\max(\padyl,\padyr)$\;
Form distance arrays $X$ and $Y$ according to \eqref{shortunsymXY}\;
Form  $X2=X^2$, $Y2=Y^2$, and $R=X2(:)+Y2$\;
\For{$r=1:\nbz$ }{
Set $z_1=z\_blocks(r)$, $z_2=z\_blocks(r+1)$\;
Calculate $grow\{1\} = \tilde{h}({11})_{pq}$, $1\le p \le nX$, $1\le q\le nY$\;
\For{$j=2:\nsy+\padyl$}{
Calculate $grow\{j\} = \tilde{h}({1j})_{p1}$, $1\le p \le nX$\;
}
Calculate : $gcol\{1\}= \tilde{h}({ij})_{11}$, $1\le i \le nX$, $1\le j\le nY$\;
\For{$q=2:\nsy+\padyr$}{
Calculate : $gcol\{q\}= \tilde{h}({i1})_{1q}$, $1\le i \le nX$\;}
\For{$j=\padyl+1:\padyl+\nsy$}{
	    Augment column of  $\Tcirc$, $gcol\{1\}$ and $grow\{j\}$, \eqref{barpaddedtoeplitz} with \eqref{Eqn:Tcircpaddedgeneral} \;}
\For{$q=\nsy+\padyr:-1:2$}{
	    Augment column of  $\Tcirc$, $gcol\{q\}$ and $grow\{1\}$, \eqref{paddedtoeplitz}  with \eqref{Eqn:Tcircpaddedgeneral} \;}
\For{$j=1:\padyl$}{
	    Augment column of  $\Tcirc$,  $gcol\{1\}$ and $grow\{j\}$, \eqref{barpaddedtoeplitz} with \eqref{Eqn:Tcircpaddedgeneral} \;}
Take FFT of  $\Tcirc$: $\hat{T}^{\mathrm{circ}}(:,:,r)=\fftt(\Tcirc)$ \;
Generate $T$ from $\Tcirc$ using \eqref{Eqn:Ttran0} and \eqref{Eqn:Ttran}\;
Take FFT of  $T$: $\hat{\tilde{T}}^{\mathrm{circ}}(:,:,r)=\fftt(T)$ \;}
$\Thatforward=\hat{T}^{\mathrm{circ}}$, $\Thattranspose=\hat{\tilde{T}}^{\mathrm{circ}}$
\end{algorithm}
\end{center}

\begin{center}
\begin{algorithm}[htb!] 
\caption{$g=\gravity(z_1,z_2, X,Y,XY,R)$\label{Alg:gravityresponse}\newline 
Entries of sensitivity matrix $G$ for the gravity problem.}
  \SetAlgoLined
  \LinesNumbered
  \KwIn{\ Depth coordinates  $z_1$ and $z_2$ for the slice\;
  $X$: Distances of $x-$coordinates from station $1$ size nx\;
  $Y$: Distances of $y-$coordinates from station $1$ size ny\;
  $XY$: the product $X(:).*Y$ which is a matrix of size $(nx+1)\times (ny+1)$\;
  $R$: the matrix of size $(nx+1)\times (ny+1)$ of entries $X(:).^{\wedge }2$ and $Y.^{\wedge }2$\;}
    \KwOut {\  Response vector $g$  of length $(nx+1) (ny+1)$\;      }
 $[nx,ny]=\textrm{size}(R)$\;
$R_1=\textrm{sqrt}(R+z_1^2)$\;
$R_2=\textrm{sqrt}(R+z_2^2)$\;
$CMX=(\textrm{log}((X(:)+R_1) ./ (X(:)+R_2))).*Y$\;
$CMY=(\textrm{log}((Y+R_1)./(Y+R_2))).*X(:)$\;
$CM5Z=\textrm{atan2}(XY,R_1z_1)z_1$\;
$CM6Z=\textrm{atan2}(XY,R_2z_2)z_2$\;
$CM56=CM5Z-CM6Z$\;
$CM=(CM56-CMY-CMX)\gamma$\;
$g=-(CM(1:nx-1,1:ny-1)-CM(1:nx-1,2:ny)-CM(2:nx,1:ny-1)+CM(2:nx,2:ny))$\;

\end{algorithm}
\end{center}

\begin{center}
\begin{algorithm}[htb!] 
\caption{$g=\mathrm{magnetic}(z_1,z_2, X,Y,R,gc))$\label{Alg:magneticresponse} \newline 
Entries of sensitivity matrix $G$ for the magnetic problem.}
  \SetAlgoLined
  \LinesNumbered
 \KwIn{\ Depth coordinates  $z_1$ and $z_2$ for the slice\;
$X$: Distances of $x-$coordinates from station  \;
$Y$: Distances of $y-$coordinates from station   \;
 $R$: Matrix   of entries $X(:)^2$ and $Y^2$\;
  $gc$ vector of constants, \cite[3]{RaoBabu:91}\;}
    \KwOut {\  Response vector $g$  of length $(\ell+1)( k+1)$\;      }
 $\ell=\textrm{length}(X)-1$;$k=\textrm{length}(Y)-1$\;
$R_1=\textrm{sqrt}(R+z_1^2)$\;
$R_2=\textrm{sqrt}(R+z_2^2)$\;
$F_1=((R_2(1:\ell,1:k)+ X(1:\ell))./(R_1(1:\ell,1:k)+ X(1:\ell))).*((R_1(2:\ell+1,1:k)+ X(2:\ell+1))./(R_2(2:\ell+1,1:k)+ X(2:\ell+1))).*((R_1(1:\ell+1,2:k+1)+ X(1:\ell))./(R_2(1:\ell+1,2:k+1)+ X(1:\ell))).*((R_2(2:\ell+1,2:k+1)+ X(2:\ell+1))./(R_1(2:\ell+1,2:k+1)+ X(2:\ell+1)))$\;     
$F_2=((R_2(1:\ell,1:k)+ Y(1:k))./(R_1(1:\ell,1:k)+ Y(1:k))).*((R_1(2:\ell+1,1:k)+ Y(1:k))./(R_2(2:\ell+1,1:k)+ Y(1:k))).*((R_1(1:\ell+1,2:k+1)+ Y(2:k+1))./(R_2(1:\ell+1,2:k+1)+ Y(2:k+1))).*((R_2(2:\ell+1,2:k+1)+ Y(2:k+1))./(R_1(2:\ell+1,2:k+1)+Y(2:k+1)))$\;
$F_3=((R_2(1:\ell,1:k)+ z2)./(R_1(1:\ell,1:k)+ z1)).*((R_1(2:\ell+1,1:k)+z1)./(R_2(2:\ell+1,1:k)+z2)).*((R_1(1:\ell+1,2:k+1)+ z1)./(R_2(1:\ell+1,2:k+1)+z2)).*((R_2(2:\ell+1,2:k+1)+ z2)./(R_1(2:\ell+1,2:k+1)+ z1))$\;
$F_4=\mathrm{atan2}(X(2:\ell+1)z_2,R_2(2:\ell+1,2:k+1).*Y(2:k+1))-\mathrm{atan2}(X(1:\ell)z_2,R_2(1:\ell+1,2:k+1).*Y(2:k+1))-\mathrm{atan2}(X(2:\ell+1)z_2,R_2(2:\ell+1,1:k).*Y(1:k))+\mathrm{atan2}(X(1:\ell)z_2,R_2(1:\ell,1:k).*Y(1:k))-\mathrm{atan2}(X(2:\ell+1)z_1,R_1(2:\ell+1,2:k+1).*Y(2:k+1))+\mathrm{atan2}(X(1:\ell)z_1,R_1(1:\ell+1,2:k+1).*Y(2:k+1))+\mathrm{atan2}(X(2:\ell+1)z_1,R_1(2:\ell+1,1:k).*Y(1:k))-\mathrm{atan2}(X(1:\ell)z_1,R_1(1:\ell,1:k).*Y(1:k))$\;
$F_5=\mathrm{atan2}(Y(2:k+1)z_2,R_2(2:\ell+1,2:k+1).*X(2:\ell+1))-\mathrm{atan2}(Y(2:k+1)z_2,R_2(1:\ell+1,2:k+1).*X(1:\ell))-\mathrm{atan2}(Y(1:k)z_2,R_2(2:\ell+1,1:k).*X(2:\ell+1))+\mathrm{atan2}(Y(1:k)z_2,R_2(1:\ell,1:k).*X(1:\ell))-\mathrm{atan2}(Y(2:k+1)z_1,R_1(2:\ell+1,2:k+1).*X(2:\ell+1))+\mathrm{atan2}(Y(2:k+1)z_1,R_1(1:\ell+1,2:k+1).*X(1:\ell))+\mathrm{atan2}(Y(1:k)z_1,R_1(2:\ell+1,1:k).*X(2:\ell+1))-\mathrm{atan2}(Y(1:k)z_1,R_1(1:\ell,1:k).*X(1:\ell))$\;           
$g=(gc(1)*\mathrm{log}(F_1)+gc(2)*\mathrm{log}(F_2)+gc(3)*\mathrm{log}(F_3)+gc(4)*F_4+gc(5)*F_5)$\;
$g=g(:)$\;

\end{algorithm}
\end{center}

\begin{center}
\begin{algorithm}[htb!]
\caption{$\bfb=\multbttb(\That, \bfx, t, \param)$\newline
This algorithm calculates the forward and transpose multiplication, $G\bfx$, or $G^T\bfx$ as described in Section~\ref{sec:circulant} using the embedding of the BTTB matrix in a BCCB matrix and the 2DFFT. The transform of \eqref{Eqn:Tcircpadded} or \eqref{Eqn:Tcircpaddedgeneral} for symBTTB and BTTB, respectively, is precomputed for the forward multiplication and provided in $\Thatforward$. The transform for \eqref{Eqn:Ttran} for the transpose, is provided in $\Thattranspose$. See Table~\ref{tab:parameters} for definitions of input parameters.\label{Alg:multbttb}}
\LinesNumbered
\KwIn{$\That$: see Table~\ref{tab:parameters}\;
$\bfx$ : vector for forward  or transpose multiplication\;
$t$ : $1$ or $2$ for forward or transpose multiplication, respectively\;
$\param$ : required parameters see Table~\ref{tab:parameters} }
\KwOut{ vector: $\bfb$ of size $m$ or $n$, for $t=1$, $2$, respectively.}
Extract parameters from $\param$ \;
Initialize zero array for $\bfb$ and $W$\;
\If{ $t==2$}{
Initialize  $W$ according to \eqref{Eqn:WTpadded} \;
Take transform of $W$:     $\hat{W}=\fftt(W)$\;
}
\For{ $j = 1:\nbz$ \% For all layers of domain}{
    \Switch{ $t$}{
        \Case{ 1}{ 
        Initialize  $W$ according to \eqref{Eqn:Wpadded}\;
        Take transform of $W$:     $\hat{W}=\fftt(W)$\;
        Form convolution \eqref{GrviaT}: $W=\real(\ifftt(\hat{T}(:,:,j)\dotstar \fftt(W)))$\;
        Extract and accumulate top left block:    $\bfb=\bfb+\reshape(W(1:\nsx,1:\nsy),m,1)$\;}
        \Case{ 2}{ 
             Form convolution \eqref{GrviaT}: $Z=\real(\ifftt(\hat{\tilde{T}}(:,:,j)\dotstar \fftt(W)))$\;
             Extract top left block and assign to output: 
            $\bfb((j-1) \nbr+1:j \nbr)=\reshape(Z(1:\nbx,1:\nby),\nbr,1)$\;}
    }
}
\end{algorithm}
\end{center}

\newpage
\bibliographystyle{plainnat}
\bibliography{/Users/rosie/Documents/Current/texfiles/cvstuff/UPRE}

\begin{thebibliography}{12}
\providecommand{\natexlab}[1]{#1}
\providecommand{\url}[1]{\texttt{#1}}
\expandafter\ifx\csname urlstyle\endcsname\relax
  \providecommand{\doi}[1]{doi: #1}\else
  \providecommand{\doi}{doi: \begingroup \urlstyle{rm}\Url}\fi

\bibitem[Boulanger and Chouteau(2001)]{BoCh:2001}
Olivier Boulanger and Michel Chouteau.
\newblock Constraints in {3D} gravity inversion.
\newblock \emph{Geophysical Prospecting}, 49\penalty0 (2):\penalty0 265--280,
  2001.
\newblock ISSN 1365-2478.
\newblock \doi{10.1046/j.1365-2478.2001.00254.x}.
\newblock URL \url{http://dx.doi.org/10.1046/j.1365-2478.2001.00254.x}.

\bibitem[Bruun and Nielsen(2007)]{bruun2007}
Christian~Eske Bruun and Trine~Brandt Nielsen.
\newblock Algorithms and software for large-scale geophysical reconstructions.
\newblock Master's thesis, Technical University of Denmark, DTU, DK-2800 Kgs.
  Lyngby, Denmark, 2007.

\bibitem[Chan and Jin(2007)]{ChanFuJin:2007}
Raymond Hon-Fu Chan and Xiao-Qing Jin.
\newblock \emph{An Introduction to Iterative Toeplitz Solvers}.
\newblock Society for Industrial and Applied Mathematics, 2007.
\newblock \doi{10.1137/1.9780898718850}.
\newblock URL \url{https://epubs.siam.org/doi/abs/10.1137/1.9780898718850}.

\bibitem[Chen and Liu(2018)]{ChenLiu:18}
Longwei Chen and Lanbo Liu.
\newblock {Fast and accurate forward modelling of gravity field using prismatic
  grids}.
\newblock \emph{Geophysical Journal International}, 216\penalty0 (2):\penalty0
  1062--1071, 11 2018.
\newblock ISSN 0956-540X.
\newblock \doi{10.1093/gji/ggy480}.
\newblock URL \url{https://doi.org/10.1093/gji/ggy480}.

\bibitem[Ha\'az(1953)]{haaz1953}
Istv\'an~B\'ela Ha\'az.
\newblock Relations between the potential of the attraction of the mass
  contained in a finite rectangular prism and its first and second derivatives.
\newblock \emph{Geophysical Transactions II}, 7:\penalty0 57--66, 1953.

\bibitem[Li et~al.(2018)Li, Chen, Chen, Dai, Zhang, Zhao, and Ling]{Li2018}
Kun Li, Long-Wei Chen, Qing-Rui Chen, Shi-Kun Dai, Qian-Jiang Zhang, Dong-Dong
  Zhao, and Jia-Xuan Ling.
\newblock Fast {3D} forward modeling of the magnetic field and gradient tensor
  on an undulated surface.
\newblock \emph{Applied Geophysics}, 15\penalty0 (3):\penalty0 500--512, Sep
  2018.
\newblock ISSN 1993-0658.
\newblock \doi{10.1007/s11770-018-0690-9}.
\newblock URL \url{https://doi.org/10.1007/s11770-018-0690-9}.

\bibitem[Pilkington(1997)]{Pilkington:97}
Mark Pilkington.
\newblock {3-D magnetic imaging using conjugate gradients}.
\newblock \emph{Geophysics}, 62\penalty0 (4):\penalty0 1132--1142, 08 1997.
\newblock ISSN 0016-8033.
\newblock \doi{10.1190/1.1444214}.
\newblock URL \url{https://doi.org/10.1190/1.1444214}.

\bibitem[Rao and Babu(1991)]{RaoBabu:91}
D.~Bhaskara Rao and N.~Ramesh Babu.
\newblock A rapid method for three-dimensional modeling of magnetic anomalies.
\newblock \emph{Geophysics}, 56\penalty0 (11):\penalty0 1729--1737, November
  1991.

\bibitem[Vogel(2002)]{Vogel:2002}
Curt Vogel.
\newblock \emph{Computational Methods for Inverse Problems}.
\newblock Society for Industrial and Applied Mathematics, Philadelphia, 2002.
\newblock \doi{10.1137/1.9780898717570}.
\newblock URL \url{http://epubs.siam.org/doi/abs/10.1137/1.9780898717570}.

\bibitem[Zhang and Wong(2015)]{ZhangWong:15}
Yile Zhang and Yau~Shu Wong.
\newblock {BTTB-based numerical schemes for three-dimensional gravity field
  inversion}.
\newblock \emph{Geophysical Journal International}, 203\penalty0 (1):\penalty0
  243--256, 08 2015.
\newblock ISSN 0956-540X.
\newblock \doi{10.1093/gji/ggv301}.
\newblock URL \url{https://doi.org/10.1093/gji/ggv301}.

\bibitem[Zhao et~al.(2018)Zhao, Chen, Chen, Liu, and Ren]{ZHAO2018294}
Guangdong Zhao, Bo~Chen, Longwei Chen, Jianxin Liu, and Zhengyong Ren.
\newblock High-accuracy {3D} {F}ourier forward modeling of gravity field based
  on the {G}auss-{FFT} technique.
\newblock \emph{Journal of Applied Geophysics}, 150:\penalty0 294 -- 303, 2018.
\newblock ISSN 0926-9851.
\newblock \doi{https://doi.org/10.1016/j.jappgeo.2018.01.002}.
\newblock URL
  \url{http://www.sciencedirect.com/science/article/pii/S0926985117301751}.

\bibitem[Zhdanov(2002)]{Zhd:2002}
Michael~S. Zhdanov.
\newblock \emph{Geophysical Inverse Theory and Regularization Problems},
  volume~36.
\newblock Elsevier, Amsterdam, 2002.

\end{thebibliography}
\end{document}